\documentclass[revision,fleqn]{elsarticle}
\usepackage[margin=2.54cm]{geometry}

\usepackage[english]{babel}
\usepackage[utf8x]{inputenc}
\usepackage[protrusion=true,expansion=true]{microtype}
\usepackage{hyperref}
\usepackage{placeins}  
\usepackage{xcolor}

\usepackage{amsmath,amssymb, empheq} 

\usepackage[displaymath, mathlines,running]{lineno}

\usepackage{algorithm}
\usepackage{algpseudocode}

\usepackage{amsfonts,nccmath,bbm,mathrsfs,mathtools}

\usepackage{etoolbox} 

\newcommand*\linenomathpatch[1]{%
  \cspreto{#1}{\linenomath}%
  \cspreto{#1*}{\linenomath}%
  \csappto{end#1}{\endlinenomath}%
  \csappto{end#1*}{\endlinenomath}%
}
\newcommand*\linenomathpatchAMS[1]{%
  \cspreto{#1}{\linenomathAMS}%
  \cspreto{#1*}{\linenomathAMS}%
  \csappto{end#1}{\endlinenomath}%
  \csappto{end#1*}{\endlinenomath}%
}

\expandafter\ifx\linenomath\linenomathWithnumbers
  \let\linenomathAMS\linenomathWithnumbers
  \patchcmd\linenomathAMS{\advance\postdisplaypenalty\linenopenalty}{}{}{}
\else
  \let\linenomathAMS\linenomathNonumbers
\fi

\linenomathpatch{equation}
\linenomathpatchAMS{gather}
\linenomathpatchAMS{multline}
\linenomathpatchAMS{align}
\linenomathpatchAMS{alignat}
\linenomathpatchAMS{flalign}

\newtheorem{thm}{Theorem}

\newdefinition{prop}[thm]{Proposition}
\newdefinition{cor}{Corollary}
\newdefinition{rmk}{Remark}
\newdefinition{define}{Definition}
\newproof{pf}{Proof}
\biboptions{sort&compress}

\definecolor{Green}{RGB}{5, 156, 16}

\definecolor{Bleu}{rgb}{0.3, 0.2, 1.0}

\definecolor{Both}{rgb}{0.55, 0.2, 0.65}

\usepackage{float}  

\newcommand{\bsym}[1]{\boldsymbol{#1}}
\newcommand{\norm}[1]{\lVert#1\rVert}

\newcommand{\Diff}{\text{Diff}(M)}

\newcommand{\SDiff}{\text{Diff}_{\mu}(M)}

\newcommand{\inprod}[2]{\langle #1,#2 \rangle}
\newcommand{\Xhi}{\mathcal{X}}
\newcommand{\bu}{\bsym{u}}
\newcommand{\id}{\text{id}}
\newcommand{\Jh}{\mathcal{J}_h}
\newcommand{\Phit}{\Phi_{\Delta t}}

\begin{document}

\begin{frontmatter}
\title{A characteristic mapping method for incompressible hydrodynamics on a rotating sphere}
\author[myaddress]{Seth Taylor\corref{correspondingauthor}}
\ead{seth.taylor@mail.mcgill.ca}

\author[myaddress]{Jean-Christophe Nave}
\address[myaddress]{\footnotesize \textit{Department of Mathematics and Statistics, McGill University, Montr\'{e}al, Qu\'{e}bec H3A 0B9, Canada}}
\cortext[correspondingauthor]{Corresponding author}

\begin{abstract}

We present a semi-Lagrangian characteristic mapping method for the incompressible Euler equations on a rotating sphere. The numerical method uses a spatio-temporal discretization of the inverse flow map generated by the Eulerian velocity as a composition of sub-interval flows formed by $C^1$ spherical spline interpolants. This approximation technique has the capacity of resolving sub-grid scales generated over time without increasing the spatial resolution of the computational grid. The numerical method is analyzed and validated using standard test cases yielding third-order accuracy in the supremum norm. Numerical experiments illustrating the unique resolution properties of the method are performed and demonstrate the ability to reproduce the forward energy cascade at sub-grid scales by upsampling the numerical solution.

\end{abstract}

\end{frontmatter}

\section{Introduction}

Direct numerical simulations of atmospheric fluid dynamics are faced with a challenge of capturing a large range of temporal and spatial scales of motion. Commonly studied fluid structures such as jets and vortices can generate multi-scale phenomena rapidly exceeding the minimum wavelength that the discretization can represent. Many of these difficulties are already manifest in the idealized dynamics governed by the barotropic vorticity equations, exhibiting key features of turbulent geophysical fluids including the advective non-linearity and dual energy cascades \cite{pedlosky1987geophysical, fjortoft1953changes}. Alternative discretization techniques capable of extending the numerical resolution of the vorticity have been proposed. \par

On a fixed Eulerian grid-based discretization, the capacity to represent fine scale features is set a priori by the spatial distribution of grid points. These techniques are thus generally more prone to numerical diffusion than their Lagrangian counterparts which allow the discretization to deform and concentrate in regions of high variation. In the context of vortex dynamics on the sphere, examples of Lagrangian methods include point vortex techniques \cite{newton2006n}, the contour dynamics and advection techniques \cite{zabusky1979contour,dritschel1988contour, dritschel1989contour}, and particle-panel methods \cite{bosler2014lagrangian}. The resulting non-uniform distribution of the discretization can however compromise the accuracy \cite{perlman1985accuracy} and ease of access to the solution throughout the entire domain. Remeshing techniques have been to developed to address this problem \cite{dukowicz1987accurate, nordmark1991rezoning, barba2005advances, magni2012accurate}, however direct interpolation of the vorticity field back onto a grid inevitably causes a numerical diffusion of fine features. In \cite{bosler2013particle, bosler2014lagrangian} an indirect remapping technique was developed which circumvents this issue by instead interpolating the inverse flow map discretized using a Lagrangian particle-panel method. The initial vorticity is then sampled at the resulting points, utilizing the conservation of absolute vorticity to define the remeshed vorticity values. The hybrid semi-Lagrangian approach of \cite{dritschel1997contour} also avoids numerical diffusion by maintaining a contour representation of the vorticity and grid-based discretization of the associated Eulerian velocity field. \par 

In this work, we present a numerical method capable of representing fine scales globally by building a spatio-temporal discretization of the inverse flow map. The method is based on the recently developed characteristic mapping (CM) method techniques for Euler's equations \cite{yin2021characteristic, yin2023characteristic} and transport on the sphere \cite{taylor2023projection}.  We utilize the advection of the vorticity, as in \cite{bosler2014lagrangian}, to compute its evolution through composition of the initial condition with the inverse map. In this form, the multiple spatial scales present in the vorticity can be attributed to those generated in the map. This property is leveraged to represent the generation of fine scale features by approximating the path formed by the inverse map as a composition of short-time sub-interval flows. The submaps are computed using the semi-Lagrangian gradient-augmented level set (GALS) method \cite{nave2010gradient}, producing a continuously differentiable approximations in a piecewise polynomial spline space through Hermite interpolation. The resulting spatio-temporal discretization formed by their composition has the capacity of representing localized oscillatory behaviour of the order $d^N$, where $d$ is the polynomial degree and $N$ is the number of compositions.  \par

Building upon the CM method techniques developed in \cite{yin2021characteristic,yin2023characteristic, taylor2023projection}, we present the design and analysis of a CM method for incompressible hydrodynamics on a rotating sphere. The numerical method and its validation serve as an essential building block for the application of this approach to other turbulent geophysical flows. A multi-grid approach is taken to compute the evolution: the inverse map is approximated using the CM method for linear advection on a spherical triangulation \cite{taylor2023projection}, coupled to a spectral method based on spherical harmonics for the vorticity via sampling. Given the ability to sample the map throughout the domain as a composition of piecewise polynomial splines, the method possesses a unique lens on the multi-scale structure of the approximated flow. We discuss and experiment with this property, demonstrating the capacity to resolve vortex structures with the expected energy cascade scalings by upsampling the solution. \par 

The paper is organized as follows:  in section \ref{sec:math_formulation} we begin by describing the mathematical formulation of the CM method and the diffeomorphism approximation technique used for the inverse map first in a general setting on smooth manifolds. Thereafter, we consider the application of the method to the incompressible Euler equations on a rotating sphere. We then give a complete description of a proposed numerical method and its algorithmic implementation in section \ref{sec:numerical_framework}. Error estimates are provided and the conservation properties are discussed. A numerical verification of the estimates is then given in section \ref{numerical_verification} where convergence tests are performed using a number of standard test cases. We conclude with numerical experiments involving a multiple zonal jet shear instability and a randomly initialized vorticity, both illustrating the unique resolution properties of the method.

\section{Mathematical Framework}
\label{sec:math_formulation}
The mapping based techniques used by the CM method are geometric in nature. Since the method discretizes a transformation of the entire domain, the numerical treatment will largely be dependent on the manifold in question.  In an effort to outline the general properties of the method which do not depend on the domain, we give a description of the mathematical framework on a compact manifold $M$. We elaborate on the use of the submap decomposition technique as a semi-discretization in time and the resolution properties of the method, phrasing some of the arbitrary resolution properties discussed in \cite{mercier2020characteristic} in terms of a relabelling symmetry held by the discretization. The section is concluded with the particular equations of motion treated in this work, describing the coupling between the map and the velocity field for the incompressible Euler equations on a rotating sphere.

\subsection{Evolution of the Inverse Flow Map}

We denote by $\varphi_{[0,t]}: M \to M$ as the forward trajectory (Lagrangian) map associated to the Eulerian fluid velocity field $\bu(t):M \to TM$. The Lagrangian velocity $\partial_t \varphi_{[0,t]}$ is related to $\bu(t)$ via
\begin{equation}
\label{eulerian_velocity}
\bsym{u} = (\partial_t\varphi_{[0,t]}) \circ \varphi_{[t,0]} \,,
\end{equation}where the map $\varphi_{[t,0]}$ is the inverse of the Lagrangian position map, i.e.
\begin{equation}
\label{inverse_map}
\varphi_{[t,0]} \circ \varphi_{[0, t]} = \varphi_{[0,t]} \circ \varphi_{[t,0]} = \text{id}_{M} \, .
\end{equation}The inverse map defines a transformation from the moving frame of reference into the Eulerian frame,  yielding the initial location $\varphi_{[t,0]}(x) = \alpha \in M$ of a fluid particle now located at $x = \varphi_{[0,t]}(\alpha) \in M$. Pointwise, the footpoint $\alpha$ is given by the solution to the ordinary differential equation equation
\begin{equation}
\label{characteristic_curves}
\dot{\bsym{\gamma}}(s) = \bsym{u}(\bsym{\gamma}(s), s) \, , \quad \bsym{\gamma}(t) = \bsym{x} \,, 
\end{equation}for the characteristic curves of $\bsym{u}$ backward in time to $s=0$. Differentiating \eqref{inverse_map} with respect to time we see that $\varphi_{[t,0]}$ satisfies an initial value problem of the form
\begin{equation}
\label{map_transport}
\partial_t\varphi_{[t,0]} + D\varphi_{[t,0]}(\bsym{u}(t)) = 0 \, ,\quad
\varphi_{[0,0]} = \text{id}_{M} \, ,
\end{equation}
in the Eulerian frame where $D\varphi_{[t,0]}\vert_{x}: T_xM \to T_{\alpha}M$ is its differential. The transport equation \eqref{map_transport} describes the advection of the fluid particle labels under the velocity \eqref{eulerian_velocity}. The paths formed by the inverse and forward maps can be decomposed into a composition of sub-interval flows. If we consider partitioning an interval of time $[0,t]$ into subdivisions $[\tau_{i}, \tau_{i+1}] \subset [0, t]$, then we have that 
\begin{equation}
\label{eq:global_map}
\varphi_{[\tau_n,0]} = \varphi_{[\tau_1, 0]} \circ \varphi_{[\tau_2, \tau_1]}  \dots \circ \varphi_{[\tau_n, \tau_{n-1}]}\,,
\end{equation}
where each of the maps $\varphi_{[\tau_{i+1}, \tau_i]}: M\to M$ are defined as the solutions to the initial value problems
\begin{equation}
\label{eq:sub_map}
\partial_t \varphi_{[t, \tau_i]}  + D\varphi_{[t,\tau_i]}(\bsym{u}(t))  = 0 \, , \quad \varphi_{[\tau_{i}, \tau_{i}]} = \id_M \, ,
\end{equation}at $t = \tau_{i+1}$. Globally, the decomposition of the solution to \eqref{map_transport} into the solutions of the manifold-valued PDEs \eqref{eq:sub_map} yields a representation of a complex deformation of $M$ as a composition of simpler and hence more accurately computed deformation maps.

\subsection{Spatio-temporal Discretization}

The approximation of the inverse map employed in \cite{yin2021characteristic, yin2023characteristic, taylor2023projection} used \eqref{eq:global_map} as a temporal discretization where each of the sub-interval flows are discretized in an interpolation space $\mathcal{V}_h \subset C^1(M,M)$. Before elaborating on the semi-Lagrangian evolution strategy used to compute each sub-interval flow and the particular choice of interpolation space, we first discuss some of the benefits gained from this approximation technique. If we denote $\mathcal{X}_{[t,s]} \in \mathcal{V}_h$ as the spatial approximation of $\varphi_{[t,s]}$, then the spatio-temporal discretization of the inverse map at $t_n$ is formed by 
\begin{equation}
\label{remapping}
\mathcal{X}_{[t_n,0]} = \mathcal{X}_{[\tau_1, 0]} \circ \mathcal{X}_{[\tau_2, \tau_1]} \circ \dots \circ \mathcal{X}_{[t_n, \tau_k]} \in \underbrace{\mathcal{V}_h\circ \mathcal{V}_h  \dots \circ  \mathcal{V}_h}_{\text{ $n_c$ times}} \,.
\end{equation}Rather than storing the map $\mathcal{X}_{[t_n,0]}$ as an element of $\mathcal{V}_h$, the decomposed maps $\mathcal{X}_{[\tau_{i+1}, \tau_{i}]} \in \mathcal{V}_h$ are stored in memory resulting in an approximation of $\mathcal{X}_{[t_n,0]}$ in an $n_c$ times composed interpolation space \eqref{remapping}. The spatial truncation of scales in the former is statically enforced by $\mathcal{V}_h$, whereas the decomposition \eqref{remapping} dynamically grows the approximation space as the flow generates finer scale features in the map. The effect of building the discretization through a composition gives this technique the capability of capturing exponentially increasing oscillatory behaviour globally without necessitating a spatial refinement of the mesh. Moreover, the degrees of freedom (d.o.f) in \eqref{remapping} only grow as $n_c \cdot |\mathcal{V}_h|$ where $|\mathcal{V}_h|$ is the d.o.f of the discretization space. This allows for a global representation of the complex multi-scale deformation generated by a turbulent fluid flow using only a linear increase in the degrees of freedom and computational resources. \par

The entire evolution is computed over a sequence of time steps $t_n$, with $n= 1, \dots, N_t$ such that the $\{\tau_j\}$ form a sub-sequence. Each submap $\mathcal{X}_{[t,\tau_j]}$ is computed over the $t_i$ such that $\tau_j \leq t_i < \tau_{j+1}$, where the remapping time $\tau_{j+1}$ can be prescribed a priori or adaptively determined. We then store $\mathcal{X}_{[\tau_{j+1}, \tau_j]} \in \mathcal{V}_h$ in memory, reinitialize $\mathcal{X}_{[\tau_{j+1}, \tau_{j+1}]} = \text{id}_{M}$, and iterate this computation.If we suppose that $M$ is Euclidean, then we can observe that the $L^{\infty}$ error for two submaps accumulates as
\begin{equation}
\label{error_accumulate}
\begin{aligned}
\norm{\varphi_{[\tau_2,0]} - \mathcal{X}_{[\tau_2,0]}}_{\infty} &\leq  \norm{(\varphi_{[\tau_1,0]} - \mathcal{X}_{[\tau_1,0]}) \circ \mathcal{X}_{[\tau_2, \tau_1]}}_{\infty}  + \norm{\varphi_{[\tau_1,0]} \circ \varphi_{[\tau_2,\tau_1]} - \varphi_{[\tau_1,0]}  \circ \mathcal{X}_{[\tau_2, \tau_1]}}_{\infty}  
\\
&\leq \norm{\varphi_{[\tau_1,0]} - \mathcal{X}_{[\tau_1,0]}}_{\infty}  + C(\tau_1)\norm{\varphi_{[\tau_2,\tau_1]} - \mathcal{X}_{[\tau_2, \tau_1]}}_{\infty}
\end{aligned}
\end{equation}
where $C(\tau_1)$ is the Lipschitz constant of $\varphi_{[\tau_1,0]}$. Since the submaps are all initialized as the identity map, which can be represented exactly, the error in $\mathcal{X}_{[t, \tau_i]}$ is monotonically increasing from zero as $t \to \tau_{i+1}$. As a result, the error accrued over the computation of the previous submaps is not carried over into the computation of the proceeding one. The approximation \eqref{remapping} thus also improves the accuracy of the method at the expense of increased memory allocation. \par

\subsection{Lie Advection}

A primary advantage of considering the inverse map as the computational quantity of interest is that pullback with $\varphi_{[t,0]}$ provides a solution operator to the homogeneous Lie advection equation 
\begin{equation}
\label{Lie_advection}
(\partial_t + \mathcal{L}_{\bu})a(t) = 0 \,, \quad  a(0) = a_0 \,, 
\end{equation}
for a general differential $k$-form $a(t) \in \Omega^{k}(M)$ where $\mathcal{L}_{\bu}: \Omega^k(M) \to \Omega^k(M)$ is the Lie derivative along the velocity field. Equation \eqref{Lie_advection} is a geometric generalization of the transport equation for more general objects than scalar fields and reduces to the advection equation in the case that $a(t) \in \Omega^0(M)$. The vorticity two-form in Euler's equations is a Lie advected quantity, regardless of the dimension.  If $a(t)$ satisfies \eqref{Lie_advection}, then we have that   

\begin{equation} 
\label{solution_op}
\frac{d}{dt} \varphi_{[0,t]}^*a(t) = \varphi_{[0,t]}^*(\partial_ta(t) + \mathcal{L}_{\bsym{u}}a(t)) = 0 \iff a(t) = \varphi_{[t,0]}^*a_0 \,,
\end{equation}
where $\varphi^*:\Omega^k(M) \to \Omega^k(M)$ denotes the pullback with the diffeomorphism $\varphi \in \Diff$. Discretizing the solution operator to \eqref{Lie_advection} rather than the solution, i.e. $\varphi_{[t,0]}^*$ rather than $a(t)$, leads to a number of advantageous properties. Since the action of sampling $a(t)$ through the inverse map can be performed arbitrarily throughout the domain via interpolation, Lie advected quantities do not necessitate an explicit discretization. Their spatial discretization can be defined instead through a chosen sampling operation, the discretization $\mathcal{X}_{[t,0]}$, and the definition of $a_0$. Since we do not interpolate the values of the transported quantity directly, this approach avoids a dissipative-type truncation error common to discretizations on a fixed Eulerian grid yet still gives ease of access to the solution throughout the entire domain since $\mathcal{X}_{[t,0]}$ is a globally-defined interpolant. The error introduced by the approximation can instead be viewed as an error in location of point evaluation and is due to the discretization respecting a relabelling symmetry. \par

\subsection{Relabelling Symmetry}

If the approximation of the map satisfies $J_{\mu}(\mathcal{X}_{[t,0]})>0$, where $J_{\mu}(\varphi)$ is the Jacobian determinant of $\varphi:M \to M$ with respect to the volume form $\mu$, then there exists a differentiable inverse to $\mathcal{X}_{[t,0]}$ by the inverse function theorem. As a result, there is a unique map $\mathcal{E}_{[t,0]} \coloneqq \mathcal{X}_{[0,t]} \circ \varphi_{[t,0]}$ which completely characterizes the error since 
\begin{equation}
\label{advective_error}
a(t) = \mathcal{E}_{[t,0]}^*\mathcal{X}_{[t,0]}^*a_0 \,,
\end{equation}
using the property $(\varphi\circ \eta)^* = \eta^*\varphi^*$ for diffeomorphisms $\varphi,\eta: M \to M$. The numerical approximation is simply a rearrangement of the initial condition, related to the true solution through the deformation $\mathcal{E}_{[t,0]}$. In turn, the method possesses a continuous form of conservation. This can be readily observed by considering the integral form of the Lie advection equation \eqref{Lie_advection}, given by the conservation law
\begin{equation}
\label{conservation_law}
\frac{d}{dt} I(a(t), S(t)) = \frac{d}{dt}\int\limits_{S(t)}a(t) = 0  
\end{equation}
where $S(t) = \varphi_{[0,t]}(S)$ and $S \subseteq M$ is a $k$-dimensional surface.  The conservation property $I(a_0, S) = I(\mathcal{X}_{[t,0]}^*a_0, \mathcal{X}_{[0,t]}(S))$ then follows due to relabelling invariance, in the sense that if $I(a(t),S(t))$ satisfies \eqref{conservation_law} then so does $I(\eta^*a(t), \eta^{-1}(S(t)))$ for any  $\eta \in \Diff$ by a change of variables. In the case of the incompressible hydrodynamics considered here, this property holds for the circulation 

\begin{equation}
\oint\limits_{\partial S} \bu^{\flat}_0 = \oint\limits_{\partial S(t)}\bu^{\flat}(t) = \int\limits_{S(t)} \bsym{\omega}(t)
\end{equation}
where $S \subset M$ is two dimensional and $\bu^{\flat}$ is the Eulerian velocity one-form. The advective nature of the error is the result of discretizing the evolution in the space of diffeomorphisms of the domain, respecting the relabelling symmetry of the conservation law \eqref{conservation_law}.
\par  

\begin{rmk}
Relabelling symmetries of the subgroup of volume-preserving diffeomorphisms $\SDiff$ play an essential role in the underlying geometry of ideal fluid flow. They are central to the correspondence between geodesics on $\SDiff$ with respect the $L^2$ metric and solutions of the incompressible Euler equations \cite{arnold1966differential}, Kelvin's circulation theorem \cite{marsden1983coadjoint}, and the Casimir invariants of two-dimensional incompressible hydrodynamics \cite{holm1998euler,cotter2013noether}. Since we discretize in the ambient space $\Diff$ as a composition of  $C^1$ piecewise polynomial splines and not directly in $\SDiff$ the numerical method does not preserve the geometric structure of ideal fluid flow. Contextualizing this discretization in relation to the geometric and structure-preserving techniques such as \cite{pavlov2011structure, gawlik2011geometric} is an interesting line of investigation however beyond the scope of this work. 
\end{rmk}

\subsection{Equations of motion}

The equations of motion which we treat in this work are those of an inviscid, incompressible fluid flow without forcing on the two-dimensional sphere $M = \mathbb{S}^2$. We note that with minor modifications, the techniques could be applied to other incompressible flows such as the single layer quasi-geostrophic equations. We consider the sphere as embedded in $\mathbb{R}^3$, centred at the origin with radius one, and rotating with constant angular velocity $\bsym{\Omega} \in \mathbb{R}^3$ with $\bsym{u}(t)$ now being the velocity field in the rotating frame. The incompressibility constraint $\text{div}(\bu) = 0$ allows for an expression of the evolution of the velocity field solely in terms of the total vorticity $\omega: \mathbb{S}^2 \to \mathbb{R}$ bringing the equations into the form of the barotropic vorticity equations
\begin{equation}
\label{barotropic_vorticity}
\partial_t \omega + \bu \cdot \nabla \omega = 0 \,, \quad \bu = -\nabla^{\perp}\Delta^{-1}(\omega - f) \,, \quad \omega(0) = \omega_0 \,,
\end{equation}
where $\nabla^{\perp}$ is the counterclockwise rotation by $\pi/2$ of the surface gradient $\nabla$ defined by the metric and $f = 2\bsym{\Omega} \cdot \bsym{n}$ is the planetary vorticity. The Lie advection of the vorticity (see  \ref{remark1} for the derivation) provides a coupling between the evolution of the velocity field to the inverse map, allowing us to express the equations of motion in the form     
\begin{subequations}
\label{eqs_of_motion}
\begin{align}
-\nabla^{\perp} \Delta^{-1}(\omega_0 \circ \varphi_{[t,0]} - f) \label{Biot-Savart} &= \bsym{u}(t)  \,,
\\
\partial_t \varphi_{[t,0]} + D\varphi_{[t,0]}(\bu(t)) &= 0 \label{map_transport2}\,,
\end{align}
\end{subequations}
We compute \eqref{Biot-Savart} using a spectral method based on spherical harmonics \cite{barrera1985vector} and \eqref{map_transport2} is computed using the CM method for linear advection on the sphere devised in \cite{taylor2023projection}.

\section{Numerical Method}
\label{sec:numerical_framework}
In this section we present the numerical methods used for the solution of the equations of motion \eqref{eqs_of_motion}. The solution algorithm involves an advect-project-reconstruct strategy coupling \eqref{Biot-Savart} to \eqref{map_transport2}. The solution to \eqref{map_transport2} is computed using the projection-based characteristic mapping method \cite{taylor2023projection} based on the semi-Lagrangian Gradient-Augmented Level Set (GALS) method \cite{nave2010gradient}. The reconstruction of the velocity field from \eqref{Biot-Savart} is performed using spherical harmonics and a projection onto the space of spherical splines. After a complete description of the solution algorithm, presented in the order of the steps taken in a single iteration, we provide error estimates serving as theoretical support for the convergence tests given in section \ref{numerical_verification}. \par


\subsection{Spatial Discretization of the Submaps}

The application of classical interpolation techniques which rely upon an underlying vector space structure are complicated by the non-linear nature of the space $C^k(M, M)$. For many manifolds of interest, embedding-based techniques \cite{grohs2019projection, grohs2013projection, gawlik2018embedding} provide a convenient treatment of this manifold-valued data approximation problem. These techniques discretize $\varphi \in C^k(M,M)$ as an embedding $M \hookrightarrow \mathbb{R}^m$ with $m > n$ and constrain the map to the manifold using a projection onto $M$. This allows for a straightforward application of well-established methods for higher-order interpolation of vector-valued functions.We consider a discretization of the submaps in a $C^1$ piecewise polynomial spline space as it yields a local and efficient evaluation of the composition \eqref{remapping} while still remaining globally differentiable. These are both desirable properties utilized by the GALS method \cite{nave2010gradient}.\par 

The spatial discretization of the inverse map in a spherical geometry devised in \cite{taylor2023projection} employed a macro-element spherical spline interpolation technique \cite{lai2007spline, alfeld1996fitting,alfeld1996bernstein}. These techniques provide a powerful computational tool to perform local Hermite interpolation on the sphere without the need to solve a linear system or construct an explicit basis. Let $\mathcal{T} = \{T_i\}_{i = 1}^{N_T}$ be a conforming spherical triangulation of the set of vertices $\mathcal{V} = \{v_i\}_{i=1}^{N_v}$ where $v_i \in \mathbb{S}^2$ for all $i = 1, \dots, N_v$.  We denote by $\mathcal{H}_d$ as the space of homogeneous trivariate polynomials and let $\mathcal{B}_d = \mathcal{H}_d\vert_{\mathbb{S}^2}$ with  $\mathcal{B}_d(\Omega) = \mathcal{B}_d\vert_{\Omega}$ such that $\Omega \subset \mathbb{S}^2$. The space of $C^r(\mathbb{S}^2)$ spherical splines of degree $d$ on the triangulation are defined by 
\begin{equation}
S^r_d(\mathcal{T}) = \left\{s \in C^r(\mathbb{S}^2) \,:\, \left.s\right\vert_{\tau} \in \mathcal{B}_d(\tau) \quad  \forall \, \tau \in \mathcal{T} \right\}\,.
\end{equation}As interpolation operator for the components of the map, we consider the spherical spline Hermite operator $\mathcal{I}_h: C^1(\mathbb{S}^2) \to S^1_2(\mathcal{T})$ on the Powell-Sabin split of $\mathcal{T}$ as detailed in \cite{lai2007spline, alfeld1996fitting, taylor2023projection}. This interpolation operator is defined with respect to only function value and gradient information at the vertices of the triangulation \cite{alfeld1996bernstein}, i.e. $\mathcal{I}_h[f] = s$ where $s \in S^1_2(\mathcal{T})$ satisfies  
\begin{equation}
\label{Hermite_problem}
s(\bsym{v}_i) = f(\bsym{v}_i) \, ,\,\,\,  D_{\bsym{e}^1_i} s(\bsym{v}_i) = D_{\bsym{e}^1_i} f(\bsym{v}_i),  \,\,\,  D_{\bsym{e}^2_{i}} s(\bsym{v}_i) = D_{\bsym{e}^2_{i}} f(\bsym{v}_i)  \,, \quad  \forall \bsym{v}_i \in \mathcal{V} \, .
\end{equation}
The interpolant is constructed as a quadratic Bernstein-B{\'e}zier polynomial in the elements of the split triangles with coefficients defined by the data \eqref{Hermite_problem}. The interpolation operator for the embedding-based spatial discretization of the inverse map on the sphere is given by
\begin{equation}
\label{Xtild}
\mathcal{J}_h : C^1(\mathbb{S}^2, \mathbb{S}^2) \to \mathcal{P}(S_2^1(\mathcal{T})^3) \subset C^1(\mathbb{S}^2, \mathbb{S}^2)\, , \quad \varphi \mapsto \mathcal{I}^{(3)}_h[\varphi]/\norm{\mathcal{I}^{(3)}_h[\varphi]} \, .
\end{equation}
Using an embedding-based approximation, the pointwise error is bound by that of $\mathcal{I}^{(3)}_h$ by a factor of 2 \cite{gawlik2018embedding}. The differential of \eqref{Xtild} also satisfies a similar error estimate, but is no longer independent of the geometry of the manifold \cite{hielscher2023approximating}, affecting the pre-asymptotic behaviour of the approximation. Nevertheless, we have that the asymptotic order of approximation inherits that of the approximation in the ambient Euclidean space. In the particular case of the quadratic spherical spline interpolation we consider here, we have that 
\begin{equation}
\label{LTE_interpolation}
\begin{aligned}
\norm{\varphi - \mathcal{J}_h[\varphi]}_{C(\mathbb{S}^2, \mathbb{S}^2)} &= \mathcal{O}(h^3) \,,
\\
\norm{D\varphi - D\mathcal{J}_h[\varphi]}_{2, \infty} &= \mathcal{O}(h^2) \,,
\end{aligned}
\end{equation}
where the norms are defined by

\begin{equation}
\label{CSS_norm}
\norm{\varphi - \mathcal{J}_h[\varphi]}_{C(\mathbb{S}^2 , \mathbb{S}^2 )} := \sup_{x \in \mathbb{S}^2 } d_{\mathbb{S}^2 }(\varphi(\bsym{x}), \mathcal{J}_h[\varphi](\bsym{x})) , \quad \norm{D\varphi}_{2,\infty} \coloneqq \sup_{\bsym{x} \in \mathbb{S}^2} \norm{D\varphi_{\bsym{x}}}_2
\end{equation}with $\norm{\,\cdot\,}_2$ being the matrix $2$-norm. Higher-order accuracy in space could be obtained using other macro-elements techniques \cite{lai2007spline, alfeld1996fitting}. The use of the Powell-Sabin interpolant however only requires data at the vertices of the triangulation while still ensuring a good degree of accuracy and global differentiability. 

\subsection{Solution Algorithm} 
\label{semi-discretization}
In this section we give a description of the steps taken during one iteration of the method, beginning with the computation of the sub-interval flow map $\mathcal{X}_{[t,\tau_k]}$. Suppose first that we know $\bsym{u}(t_i)$ for $0\leq i \leq n$,  we begin by extrapolating the velocity field in time as  
\begin{equation}
\label{velocity_lagrange}
\tilde{\bsym{u}}(\bsym{x},t) = \sum_{i= 0}^p \ell_i(t) \bsym{u}^{n-i}(\bsym{x})
\end{equation}where $\ell_i(t)$ are Lagrange basis functions. The number $p$ is chosen such that this extrapolation is of the same order of approximation in time as the trajectory computations. Using this approximation of the velocity field, we then perform one step of the backward semi-Lagrangian Gradient-Augmented Level Set (GALS) \cite{nave2010gradient} method to compute the evolution of the submap. Each iteration updates the submap as
 
\begin{subequations}
\label{CMmethod}
\begin{align}
\mathcal{X}_{[t_{n+1}, t_n]}(\bsym{x}) &= \Phi_{\Delta t}(\tilde{\bsym{u}},\bsym{x}) \label{numerical_int}\, ,
\\
\mathcal{X}_{[t_{n+1}, \tau_k]}(\bsym{x}) &= \mathcal{J}_h [\mathcal{X}_{[t_n, \tau_k]} \circ \mathcal{X}_{[t_{n+1}, t_n]}](\bsym{x}) \label{interpolation_step}\, ,
\end{align}
\end{subequations}
where $\Phi_{\Delta t} : \mathfrak{X}(\mathbb{S}^2 \times \mathbb{R}) \times \mathbb{S}^2 \to \mathbb{S}^2$ is a numerical integration scheme used to perform the trajectory computations and $\mathcal{X}_{[\tau_k,\tau_k]}$ is initialized as the identity map. 

\subsubsection{Trajectory Computations}

The trajectory computations \eqref{numerical_int} providing the approximated values of the submap are performed in Cartesian coordinates using the classical RK4 numerical integration scheme applied to the ordinary differential equation

\begin{equation}
\label{characteristics}
\dot{\bsym{\gamma}}(t)= \tilde{\bsym{u}}(\bsym{\gamma}(t),t)\, , \quad \bsym{\gamma}(t_{n+1}) = \bsym{x} \, ,
\end{equation}backwards in time to $t_n$. The intermediate stages of the integration scheme are projected back on the sphere to evaluate the velocity field. Since the radial distance of the trajectory from the sphere is on the same order as the local error of the integration scheme, the accuracy of the trajectory computation will not be compromised through this added projection step \cite{hairer2006structure}. At each iteration the foot points of a four point $\epsilon$-difference stencil about each vertex are computed, providing approximate interpolation data to the interpolation operator \eqref{Hermite_problem} for the components of the map. These are introduced by pre-computing the positions of four stencil points $\varepsilon^i_{\pm,\pm} \in \mathbb{R}^3$ at each vertex $v_i$ as $\varepsilon^i_{\pm,\pm} = \pi_{v_i}^{-1}(\bsym{v}_i \pm \epsilon\bsym{\gamma}_{v_i}^1 \pm \epsilon\bsym{\gamma}_{v_i}^2, \bsym{v}_i \pm \epsilon\bsym{\gamma}_{v_i}^2 \pm \epsilon\bsym{\gamma}_{v_i}^1)$ where $\pi^{-1}_{v_i}$ is the inverse of the tangent plane projection at $v_i$ and $\{\bsym{\gamma}_{v_i}^1, \bsym{\gamma}_{v_i}^2\}$ form a local orthonormal basis for the tangent plane at the vertex. Applying the numerical integration scheme from these initial positions, we obtain the points

\begin{equation}
\label{foot_points}
\mathcal{X}_{[t_{n+1}, t_n]}(\varepsilon^i_{\pm, \pm}) \coloneqq \bsym{x}^{n+1,i}_{\pm,\pm} =\Phi_{\Delta t}(\tilde{\bsym{u}}, \varepsilon^i_{\pm, \pm}) \in \mathbb{R}^3 \,, \quad \forall \bsym{v}_i \in \mathcal{V}.
\end{equation}

\subsubsection{Approximated Hermite Data}

The interpolation step \eqref{interpolation_step} is facilitated by forming an approximation of the interpolation information \eqref{Hermite_problem} in each component of the map. After integrating for the footpoints \eqref{foot_points} we evaluate the map from the previous location at these points and approximate the Hermite data using the following $\epsilon$-finite-difference stencils 
\begin{equation}
\begin{aligned}
\mathcal{X}_{[t_{n+1},\tau_k]}(\bsym{v}_i) & \approx \frac{1}{4}\left( \mathcal{X}_{[t_{n},\tau_k]}(\bsym{x}^{n+1,i}_{-,-}) + \mathcal{X}_{[t_{n},\tau_k]}(\bsym{x}^{n+1,i}_{+,-}) +  \mathcal{X}_{[t_{n},\tau_k]}(\bsym{x}^{n+1,i}_{-,+}) + \mathcal{X}_{[t_{n},\tau_k]}(\bsym{x}^{n+1,i}_{+,+})\right) \,, 
\\
\bsym{\gamma}_{v_i}^1 \cdot \nabla \mathcal{X}_{[t_{n+1},\tau_k]}(\bsym{v}_i) &\approx \frac{1}{4\epsilon}\left( \mathcal{X}_{[t_{n},\tau_k]}(\bsym{x}^{n+1,i}_{+,-}) - \mathcal{X}_{[t_{n},\tau_k]}(\bsym{x}^{n+1,i}_{-,-}) +  \mathcal{X}_{[t_{n},\tau_k]}(\bsym{x}^{n+1,i}_{+,+}) - \mathcal{X}_{[t_{n},\tau_k]}(\bsym{x}^{n+1,i}_{-,+})\right) \,,
\\
\bsym{\gamma}_{v_i}^2 \cdot \nabla \mathcal{X}_{[t_{n+1},\tau_k]}(\bsym{v}_i) & \approx \frac{1}{4\epsilon}\left( \mathcal{X}_{[t_{n},\tau_k]}(\bsym{x}^{n+1,i}_{-,+}) - \mathcal{X}_{[t_{n},\tau_k]}(\bsym{x}^{n+1,i}_{-,-}) +  \mathcal{X}_{[t_{n},\tau_k]}(\bsym{x}^{n+1,i}_{+,+}) - \mathcal{X}_{[t_{n},\tau_k]}(\bsym{x}^{n+1,i}_{+,-})\right) \,.
\end{aligned}
\end{equation}
Choosing $\varepsilon$ sufficiently small, these approximations will not compromise the local truncation error of the interpolant \cite{taylor2023projection}. We refer to \cite{mercier2020characteristic, nave2010gradient} for further details on the GALS method and \cite{taylor2023projection} for its implementation for spherical splines. Finally, a new interpolant is formed as in \eqref{interpolation_step} by projecting back onto the space $\mathcal{P}(S^1_2(\mathcal{T}_{PS})^3)$ using $\mathcal{J}_h$, defining the next submap.\par

\subsubsection{Velocity Field Reconstruction}
\label{sec:velocity_field}
The remaining part of the solution strategy involves a reconstruction of the velocity field at $t_{n+1}$ through a sampling of  $\omega_0 \circ \mathcal{X}_{[t_{n+1},0]}$. The spherical triangulation $\mathcal{T}$ used for the discretization of the inverse flow map does not constrain the discretization of the velocity field. The solution strategy therefore involves a two-grid approach separating the dynamics driving the evolution of the map and the deformation of the domain which it induces. This consequently gives the reconstruction step flexibility for the way in which the Poisson solve is performed. In this work we consider a reconstruction based on a spherical harmonic expansion of the vorticity. \par 

 The spatial discretization of the velocity field is formed by a spherical harmonic expansion of the stream function. An exact band-limited spherical harmonic transform was devised by McEwen and Wiaux \cite{mcewen2011novel} based on a sampling at the grid points

\begin{equation}
\label{grid_points}
(\lambda_q, \theta_p) = \left(\frac{2\pi q}{2L-1}, \frac{\pi p}{L-1} \right)
\end{equation}where $q = 0,1,\dots, 2L-2 $ and $p = 0,1, \dots, L-1$ and $L$ is the band-limit of the spherical harmonic representation. These points define the dynamics grid $U_L$. Using the current submap $\mathcal{X}_{[t,\tau_k]}$ and the submaps $\mathcal{X}_{[\tau_k,\tau_{k-1}]} \dots \mathcal{X}_{[\tau_1,0]}$ stored in memory, we reconstruct the velocity field by first transporting the absolute vorticity scalar field on these points as 
\begin{equation}
\label{rel_vort_np1}
\zeta^{n+1}(\bsym{x}) = (\zeta_0 + f) \circ \mathcal{X}_{[t_{n+1},0]}(\bsym{x}) - f(\bsym{x}) \,,
\end{equation}
where $\mathcal{X}_{[t_{n+1},0]}$ is formed as in \eqref{remapping}. Based on these samples, we expand the vorticity in a spherical harmonic basis using the transform devised in \cite{mcewen2011novel}. Since the spherical harmonics are eigenfunctions of the spherical Laplacian with eigenvalues $-\ell(\ell + 1)$ we obtain the stream function as
\begin{equation}
\label{psi_np1}
\psi^{n+1} = \sum_{\ell = 1}^{L-1} \sum_{m= -\ell}^{\ell} \frac{\hat{\zeta}^{n+1}_{\ell,m}}{\ell(\ell + 1)} Y^{m}_{\ell} \, ,
\end{equation}
where the coefficients $\hat{\zeta}^{n+1}_{\ell,m}$ are the spherical harmonic coefficients of the relative vorticity \eqref{rel_vort_np1}.  We then differentiate the basis functions directly to recover the velocity field using
\begin{equation}
\label{u_ang_mom}
\bsym{u} = -\nabla^{\perp} \psi  = \nabla \psi \times \bsym{x} = -i(\bsym{L}_x + \bsym{L}_y + \bsym{L}_z )\psi\, , \end{equation}
where the angular momentum operators act on the spherical harmonic basis functions as defined in \ref{sph_harmonic_appendix}.  This representation of the velocity field does not suffer from any coordinate singularities induced by the spherical coordinate representation of the stream function, since the operators \eqref{L_operator} are well-defined on spherical harmonics. The construction can be shown to be equivalent to an expansion of  the velocity field in terms of the divergence-free vector spherical harmonic \cite{barrera1985vector}.\par 

The stream function \eqref{psi_np1} is then supplied to \eqref{u_ang_mom}, yielding a vectorial spherical harmonic representation of the velocity field. We then project the components of the velocity field onto $S^1_2(\mathcal{T}^{u}_{PS})^3$ where $\mathcal{T}^{u}$ is a spherical triangulation of the grid points \eqref{grid_points}. Altogether, the reconstruction of the velocity field at time $t_n$ can be written as
\begin{equation}
\label{u_sph}
\bsym{u}^{n+1} = \mathcal{I}^{(3)}_h\left[-i\bsym{L}\psi^{n+1}\right] \, ,
\end{equation}
where $\psi^{n+1}$ is given by \eqref{psi_np1}. The derivative values for the components of the velocity field needed to perform the projection \eqref{u_sph} onto $S^1_2(\mathcal{T}^u_{PS})^3$ are computed from the spherical harmonic coefficients using a rotation of the angular momentum operator. 

\begin{rmk}
Note that the divergence-free constraint on the velocity is not enforced away from the vertices. This condition could be enforced directly by first interpolating the stream function and then taking the rotated gradient, although this would yield a $C^0$ approximation of the velocity field. We have chosen to only approximate the divergence-free condition in order to retain the $C^1$ regularity. \par
\end{rmk}

\subsubsection{Summary of Implementation}
\label{sec:submap}

The implementation of the solution algorithm incorporating the submap decomposition can be written in the following pseudo-code format \ref{euler-evolution}. We refer to $\text{CM-Submap}$ as the application of a semi-Lagrangian transport step as used in \cite{taylor2023projection} and  the choice of when to perform a submap decomposition step is written generically as some Boolean criterion $\mathcal{C}: A \to \{0,1\}$ over a parameter space $A$. In \cite{yin2021characteristic,yin2023characteristic} this criterion was provided by an error tolerance on the Jacobian determinant of the submap, penalizing the deviation of the approximation from the volume-preserving diffeomorphism group.  Here, we opt for a statically enforced remapping criterion, specifying the number of steps before remapping a priori. The possibilities for adaptivity with the technique of submap decomposition will be investigated in our future work.      
\begin{algorithm}
\label{euler-evolution}
\caption{Evolution Algorithm}
\textbf{Input}: Initial vorticity $\omega_0$, final time $T$, time step $\Delta t$, remapping criterion $C$.
\\
\textbf{Output}: List  $[\mathcal{X}_{[\tau_1, 0]}, \mathcal{X}_{[\tau_2, \tau_1]}, \dots \mathcal{X}_{[T, \tau_n]}]$ 
\\
Initialization: submaps = [  ], $(U = [\bsym{u}_0, \bsym{u}^{\Delta t}, \bsym{u}^{2\Delta t}], \mathcal{X}_{[2\Delta t, 0]})$, $t \gets 2\Delta t$.
  \begin{algorithmic}[1]
  \While{$t < T$}
      \State $\tilde{\bsym{u}} \gets [\bsym{u}^{t-2\Delta t}, \bsym{u}^{t-\Delta t}, \bsym{u}^{t}]$ \Comment{using \eqref{velocity_lagrange}}

      \State $\mathcal{X}_{[t + \Delta t, t_i]}$ $\gets$ CM-Submap($\tilde{\bsym{u}}$, $\mathcal{X}_{[t, t_i]}$) \Comment{transport inverse map}
      \State $\bsym{u}^{t+\Delta t} = \mathcal{I}_h^3[\bsym{L}\Delta^{-1}(\omega_0 \circ \mathcal{X}_{[t + \Delta t,0]} - f)]$ \Comment{using \eqref{remapping} and \eqref{u_sph}}
      \State $U = [\bsym{u}^{t-\Delta t}, \bsym{u}^{t}, \bsym{u}^{t + \Delta t}]$ \Comment{update velocity field list}
      \State $t \leftarrow t + \Delta t$
      \If{$C(A) = 1$}
          \State submaps $\gets \text{submaps} \cup \mathcal{X}_{[t,t_i]}$,  $t_i \gets t$, $\mathcal{X}_{[t_i, t_i]} \gets \text{id}_M$
      \EndIf
  \EndWhile
  \end{algorithmic}
\end{algorithm}

\subsection{Error Estimates}
\label{error_estimates}

In this section we provide error estimates for the method including an analysis on the effect of the decomposition \eqref{remapping} for a fixed remapping strategy. These serve as theoretical justification for the accuracy observed in the forthcoming convergence tests. To this end, we first consider the approximation of the inverse map resulting from a modified equation for the velocity field $\tilde{\bsym{u}}$ defined by \eqref{velocity_lagrange}. This defines an evolution equation for the map $\tilde{\varphi}_{[t,0]}$ computed as the solution to the following initial value problem
\begin{equation}
\label{modified_advection}
\begin{aligned}
\partial_t\tilde{\varphi}_{[t,0]} + D\tilde{\varphi}_{[t,0]}(\tilde{\bsym{u}}(t)) &= 0\, ,\\
\tilde{\varphi}_{[0,0]} &= \text{id}_{\mathbb{S}^2} \, .
\end{aligned}
\end{equation}Let $t_n = n \Delta t$ be time steps forming an uniform partition of $[0,T]$ and suppose we use a fixed remapping strategy, where a submap is computed and stored at each $\tau_k = k\Delta\tau$ where $\Delta \tau = m\Delta t$ for some whole number $m > 1$.  We suppress the dependence on $\tilde{\bsym{u}}$ and let $\Phi_{\Delta t}:\mathbb{S}^2\times \mathbb{R} \to \mathbb{S}^2$ be the map approximating the departure points from the modified velocity field over one time step, which lets us write $\mathcal{X}_{[t_1,0]} = \mathcal{J}_h[\Phi_{\Delta t}(t_1)]$.  We can decompose the error over one subinterval $[\tau_k, t_n]$ as 

\begin{equation}
\begin{aligned}
\norm{\tilde{\varphi}_{[t_n, \tau_k]} - \mathcal{X}_{[t_n,\tau_k]}}_{\infty} &= \norm{\tilde{\varphi}_{[t_n,\tau_k]} - \mathcal{J}_h[(\mathcal{X}_{[t_{n-1},\tau_k]}-\tilde{\varphi}_{[t_{n-1},\tau_k]}) \circ \Phit(t_n)] - \mathcal{J}_h[\tilde{\varphi}_{[t_{n-1},\tau_k]}\circ \Phi_{\Delta t}(t_n)]}_{\infty} 
\\
& \leq \norm{\tilde{\varphi}_{[t_n,\tau_k]} - \Jh[\tilde{\varphi}_{[t_{n-1},\tau_k]} \circ \Phit(t_n)]}_{\infty} + \norm{\Jh}\norm{\mathcal{X}_{[t_{n-1},\tau_k]}- \varphi_{[t_{n-1},\tau_k]}}_{\infty}
\end{aligned}
\end{equation}
The first term can then be decomposed further as
\begin{equation}
\label{second_estimate}
\begin{aligned}
\norm{\tilde{\varphi}_{[t_n,\tau_k]} - \Jh[\tilde{\varphi}_{[t_{n-1},\tau_k]} \circ \Phit(t_n)]}_{\infty} &\leq \norm{\tilde{\varphi}_{[t_{n-1},\tau_k]}\circ \Phit(t_n) - \Jh[\tilde{\varphi}_{[t_{n-1},\tau_k]} \circ \Phit(t_n)]}_{\infty} \\
&+ C_{n-1,k}\norm{\tilde{\varphi}_{[t_n,t_{n-1}]} - \Phit(t_n)}_{\infty} \,.
\end{aligned}
\end{equation}
where $C_{n-1,k}$ is the Lipschitz constant of $\tilde{\varphi}_{[t_{n-1},\tau_k]}$. Using a s-stage RK integration scheme for the departure point computations the second term in \eqref{second_estimate} is $\mathcal{O}(\Delta t^{s+1})$ and using quadratic spherical spline interpolation the first term is $\mathcal{O}(C(t_n)(h^3 +  \Delta t h^2))$ \cite{taylor2023projection} with the constant $C(t_n) \sim \norm{\tilde{\varphi}_{[t_n,\tau_k]}}_{3,\infty}$ where  $\norm{\cdot}_{3,\infty}$ is the Sobolev $W^{3,\infty}$ norm \cite{lai2007spline}. The global error can be decomposed as
\begin{equation}
\label{remapping_error}
\begin{aligned}
\norm{\tilde{\varphi}_{[t_n,0]} - \mathcal{X}_{[t_n,0]}}_{\infty} &= \norm{\tilde{\varphi}_{[t_n,0]}- \mathcal{X}_{[\tau_1,0]}\circ \mathcal{X}_{[\tau_2, \tau_1]}\cdot \dots \circ \mathcal{X}_{[t_n, \tau_k]}}_{\infty} 
\\
&\leq  \norm{(\tilde{\varphi}_{[t_{k},0]} - \mathcal{X}_{[\tau_k,0]}) \circ \mathcal{X}_{[t_{n}, \tau_k]}}_{\infty} + C_k\norm{\tilde{\varphi}_{[t_n,\tau_k]} - \mathcal{X}_{[t_{n}, \tau_k]}}_{\infty}  
\\
&\leq \norm{\tilde{\varphi}_{[\tau_k,0]} - \mathcal{X}_{[\tau_k,0]}}_{\infty} + C_k\norm{\tilde{\varphi}_{[t_n,\tau_k]} - \mathcal{X}_{[t_{n}, \tau_k]}}_{\infty}   \,.
\end{aligned}
\end{equation}
Letting $t_n = \tau_{k+1}$ and considering the errors accrued over each sub-interval $[\tau_k, \tau_{k+1}]$ we get the following bound
\begin{equation}
\begin{aligned}
\norm{\tilde{\varphi}_{[t_n, 0]}- \mathcal{X}_{[t_n,0]}}_{\infty} &\leq \sum_{i=1}^k C_i\norm{\tilde{\varphi}_{[\tau_{i+1}, \tau_i]} - \mathcal{X}_{[\tau_{i+1}, \tau_i]}}_{\infty} 
\\
&\lesssim t_n\left(\max_{i = 1,\dots ,k} \norm{\tilde{\varphi}_{[\tau_{i},\tau_{i-1}]}}_{3,\infty} \cdot (h^2 + h^3/\Delta t) + \Delta t^s\right)\,,
\end{aligned}
\end{equation}
 Since the $\tilde{\varphi}_{[\tau_{k+1},\tau_{k}]}$ start from the identity map at $\tau_k$ we can say that $\norm{\tilde{\varphi}_{[\tau_{k+1},\tau_{k}]}}_{3,\infty} = \mathcal{O}(\Delta \tau)$ which gives an extra parameter to $h$ and $\Delta t$ to control the error of the approximation.

\par 
Incorporating the approximation of the velocity field, we introduce a smoothing error due to an undersampling in \eqref{rel_vort_np1} for the implementation of the Biot-Savart law \cite{yin2021characteristic}. Let $\mathcal{Q}_{L}: L^2(\mathbb{S}^2) \to L^2(\mathbb{S}^2)$ be the projection operator onto the first $\ell \leq L$ spherical harmonics. Before the projection onto the space of spherical splines we have that $\tilde{\bsym{u}} = \mathcal{Q}_L[\tilde{\bsym{u}}]$ for $L$ larger than the band-limit defining the sampling grid \eqref{grid_points}. This allows us to split the error introduced during the reconstruction step as follows 
\begin{equation}
\label{velocity_error_separation}
\norm{\bsym{u} - \tilde{\bsym{u}}}_{\infty} \leq \norm{\mathcal{Q}_L[\bsym{u}] - \tilde{\bsym{u}}}_{\infty} + C\norm{(I - \mathcal{Q}_L)[\bsym{u}]}_{H^s} \,,
\end{equation}
due to the Sobolev embedding $H^s \hookrightarrow C^0$ for $s > 1$ since we are working in two-dimensions. We can control the first term in \eqref{velocity_error_separation} with the approximation of the map, whereas we must assume that $L$ is taken large enough such that the contribution of the second term is negligible in comparison. This assumes that the analytic velocity field has sufficient decay in its energy spectrum and is justifiable for two-dimensional turbulence \cite{boffetta2012two}. Assuming that we have taken $L$ large enough such that this is the case, the global accuracy of the numerical method is described by the following theorem.  

\begin{thm}
\label{theorem1}
Let $\tilde{\bsym{u}}$ be defined by \eqref{velocity_lagrange} using a $p^{th}$ order Lagrange interpolant in time. Using an $s$-stage RK integration scheme for the departure point computations with time steps of size $\Delta t$ and using remapping steps of size $\Delta \tau$, the global error for the inverse map to final integration time $T$ is given by

\begin{equation}
\label{CM_method_error}
\norm{\mathcal{X}_{[T,0]} - \varphi_{[T,0]}}_{C^{0,\alpha}} = \mathcal{O}(T\Delta t^s + T\Delta \tau\min(h^3\Delta t^{-1}, h^2) + T\Delta t^p) \,.
\end{equation}
\end{thm}
The proof is included  in  \ref{convergence_proof} for the sake of completeness and follows from standard estimates for semi-Lagrangian schemes and elliptic regularity. The accuracy of the method is assessed using the supremum norm which is expected to have the same order of accuracy as  \eqref{CM_method_error} taking $\alpha$ arbitrarily small. Note then that the error in the map \eqref{CM_method_error} bounds the vorticity error in the supremum norm with a multiplicative factor given by the Lipschitz constant of the vorticity initial condition. \par 
In our numerical results we also assess the convergence of the method with respect to the energy and enstrophy. We note however that, due to the unique geometric structure of two-dimensional incompressible fluid flow, there exists an infinite number of other conserved quantities \cite{arnold2021topological}. In particular, by a change of variables with the volume-preserving map $\varphi_{[t,0]}$, for any measurable function $h: \mathbb{R} \to \mathbb{R}$ the integrals

\begin{equation}
\label{casimir_invariants}
I_h[\omega(t)] = \int_{\mathbb{S}^2} h(\omega(t))\mu\,,
\end{equation}are constant in time. If $h$ is Lipschitz continuous, then it can be shown that the conservation error for the invariants \eqref{casimir_invariants} satisfies the following bound

\begin{equation}
|I_h[\tilde{\omega}] - I_{h}[\omega_0]| \leq C \min \{\norm{J_{\mu}(\mathcal{X}_{[t,0]}) -1}_{\infty}, \norm{\mathcal{X}_{[t,0]} - \varphi_{[t,0]}}_{\infty}\} \,,
\end{equation}and is thus controlled by the error in the map. This holds similarly for the conservation of energy, which we do not expect to be conserved numerically since time-reversibility is not enforced. However, the energy is still controlled by the error for the inverse flow map since  

\begin{equation}
\begin{aligned}
|\norm{\tilde{\bsym{u}}(t)}^2_{L^2} - \norm{\bsym{u}(0)}^2_{L^2}| 
&= |(\tilde{\omega}(t), \Delta^{-1}\tilde{\omega}(t)) -  (\omega(t), \Delta^{-1}\omega(t))|
\\
&\leq \norm{\Delta^{-1}\tilde{\omega}(t)}_{L^2}\norm{\tilde{\omega}(t) - \omega(t)}_{L^2} + \norm{\omega(t)}_{L^2} \norm{\Delta^{-1}(\omega(t) - \tilde{\omega}(t))}_{L^2} 
\\
&\leq C \norm{\mathcal{X}_{[t,0]} - \varphi_{[t,0]}}_{C^0} \,.
\end{aligned}
\end{equation}

\section{Numerical Verification}
\label{numerical_verification}
In this section we provide a verification of our implementation using convergence tests of the method. The error estimates given in section \ref{error_estimates} are affirmed using test cases defined by initial vorticities consisting of a Rossby-Haurwitz wave, a Gaussian vortex, and a steady zonal jet. 

\subsection{Implementation Details}
\label{additional_details}
The numerical tests were implemented in Python and run on a Linux workstation with an Intel core i5-8250U (8 logical processors) with 16 GB of RAM. The spherical triangulations were constructed using the Python package Stripy \cite{moresi2019stripy} which provides a wrapper to the package STRIPACK \cite{renka1997algorithm}. The point in triangle querying was performed using the Python binding to the lib-igl package \cite{libigl}. The spherical harmonic transformations were performed using the Python binding to the SSHT package \cite{mcewen2011novel}. The inverse map is discretized using an icosahedral discretization of the sphere (see table \ref{ico_table}). We note however that the formulation and implementation are essentially agnostic to the particular spherical triangulation of the backward characteristic map. \par

\begin{table}[!h]
\centering
\begin{tabular}{|c||c|c|c|c|c|c|c|c|c|}
\hline
$k$ & 0 & 1 & 2 & 3 & 4 & 5 & 6 & 7 & 8 \\
\hline\hline 
$N_v$  & 12 & 42 & 162 & 642 & 2562 & 10242 & 40062 & 163842 & 655362\\
\hline
$N_{\Delta}$ & 20 & 80 & 320 & 1280 & 5120 & 20480 & 81920 & 327680 & 1310720 \\
\hline
$h$ & 1.10715 & 0.62832 & 0.32637 & 0.16483 & 0.08263 & 0.04134 & 0.02067 & 0.01034 & 0.00517 \\ \hline 
\end{tabular}
\caption{Number of vertices ($N_v$), simplices ($N_{\Delta}$), and maximum edge length $h$ for the $k^{\text{th}}$ refinement of the icosahedral discretization of the sphere.}
\label{ico_table}
\end{table}

The computational cost of the algorithm can be broken down into three primary components: 1) the transport of the map, 2) the evolution of the vorticity, and 3) the computation of the velocity field. Since we use an explicit time-stepping, the transport of the submap is efficient with the dominant contribution to the computational time at each iteration attributed to the evaluation of $\mathcal{X}_{[t_n,\tau_k]}$ at the footpoints of the $\epsilon$-difference stencils. As this requires the point in triangle querying it results in a $\mathcal{O}(N_{v})$ operation without any additional data structures on the triangulation. The evolution of the vorticity field at $N$ points therefore requires an $\mathcal{O}(N_c \cdot N)$ operation.  In turn, the computation of the velocity field requires an $\mathcal{O}(N_c \cdot L^2)$ operation in the sampling of the vorticity along with a $\mathcal{O}(L^3)$ operation for the spherical harmonic transform using the SSHT package \cite{mcewen2011novel} and is the dominant contribution to the computational time. Overall, the computational cost associated to the algorithm is of the order $\mathcal{O}(N_t(N_c L^2 + L^3 + N_v))$. \par 
Our implementation was written in a high-level language and the tests were performed on a laptop computer for the purpose of verifying convergence of the method. A performance optimization of the algorithm is beyond the scope of this work, however run times are modest for our purposes. Wall-clock times for the last two data points in \eqref{fig:euler_test_all} were respectively measured to be approximately $17$ and $145$ minutes without remapping and $24$ and $308$ minutes with remapping. Computational time associated to each of the operations involved in the implementation could be reduced by incorporating tree data structures on the triangulation, parallelization evaluating the map and for the foot point calculations, along with a lower-level implementation. A performance optimization was devised by utilizing the uniformity of the velocity field grid for a faster containing triangle querying strategy during the foot point calculations. The spherical triangulation resulting from the vertices \eqref{grid_points} is separated into cells $C_{i,j} \coloneqq \{(\lambda_{i}, \theta_{j}), (\lambda_{i+1}, \theta_{j}), (\lambda_{i}, \theta_{j+1}), (\lambda_{i+1}, \theta_{j+1})\}$ where $\lambda_i = \lambda_{i+1}$ if $j = 0,L$. The $C_{i,j}$ are then split from $(\lambda_i, \theta_j)$ to $(\lambda_{i+1}, \theta_{j+1})$ along a great circle arc, yielding two triangles within each cell. We omit the cells with a vertex at either pole and an array $G$  of size $[L,2L, 2]$ is then defined such that $G[j,i,0 \,(1)]$ yields the list of the vertices of the bottom (top) triangle within $C_{i,j}$. Let $(\lambda_q, \theta_q)$ be a query point on this triangulated mesh for the velocity field. The row index of the containing triangle is given simply by $i_q = \lfloor \lambda_q/\Delta x\rfloor$  and a preliminary column index is given by $j_q = \lfloor \theta_q/\Delta x \rfloor$ where $\Delta x = 2\pi/L$. The column index does not in general give the exact containing triangle since the great circle arcs connecting adjacent vertices of different longitude do not transform into straight lines in the $(\lambda, \theta)$ parametric space. The containing triangle can however be determined using the position within the cell $C_{i_q, j_q}$ and the sign of the distance from the plane containing the nearest great circle arc along the diagonal or along the two top edges of the cell. \par

\subsection{Error Norms}
Based on the error estimates provided in section \ref{error_estimates}, we assess the accuracy of the method using an approximation of the following error norms:

\begin{equation}
\label{error_norms}
\begin{aligned}
\text{Vorticity error} &\coloneqq  \frac{\norm{\omega_0 \circ \mathcal{X}_{[T,0]} - \omega (\cdot,T)}_{L^{\infty}(\mathbb{S}^2)}}{\norm{\omega (\cdot,T)}_{L^{\infty}(\mathbb{S}^2)}}\, ,
\\
\text{Enstrophy conservation error} &\coloneqq \frac{\norm{\omega_0 \circ \mathcal{X}_{[T,0]}}^2_{L^2(\mathbb{S}^2)} - \norm{\omega_0}^2_{L^2(\mathbb{S}^2)}}{\norm{\omega_0}^2_{L^2(\mathbb{S}^2)}} \, ,
\\
\text{Energy conservation error} & \coloneqq \frac{\norm{\bsym{u}^n}^2_{L^2(\mathbb{S}^2)} - \norm{\bsym{u}_0}^2_{L^2(\mathbb{S}^2)}}{\norm{\bsym{u}_0}^2_{L^2(\mathbb{S}^2)}} \,.
\end{aligned}
\end{equation}
The sup-norm errors are approximated using a sampling of the grid points \eqref{grid_points} for a band-limit of $L= 1000$. The $L^2(\mathbb{S}^2)$ norm is approximated using the spherical harmonic coefficients of the vorticity. The kinetic energy at time $t = t_n$ is approximated as 

\begin{equation}
\norm{\bsym{u}^n}^2_{L^2(\mathbb{S}^2)} = (\nabla^{\perp}\psi^n,\nabla^{\perp}\psi^n)_{\mathbb{S}^2} = (\omega^n,\psi^n)_{\mathbb{S}^2} = \sum_{\ell=1}^L\sum_{m = -\ell}^{\ell} \frac{|\hat{\omega}^n_{\ell,m}|^2}{\ell(\ell+1)} \,.
\end{equation}
We note that it is common to assess the accuracy of the method in the $\ell^{\infty}$ and $\ell^2$ norms. Here we have chosen a finer approximation of the continuous error measures in an effort to emphasize the functional definition of the inverse map and the vorticity.

\subsection{Convergence Tests}

 We demonstrate the accuracy of the method, described by \eqref{CM_method_error}, by refining $T/N_t = \Delta t$ proportionally to $h$ and $L$. In particular, we consider $N_t = 2^{k+2}$ and $L = 2^{k+3}$ where $k$ is the number of refinements of the icosahedral discretization ranging from 1-6.  A value of $\varepsilon = 10^{-5}$ is chosen for the $\varepsilon$-difference stencils in each test which effectively limits the machine precision to approximately $10^{-12}$. We perform the convergence tests both with and without submap decomposition. The time steps used for remapping were taken to be $\Delta \tau = 10 \Delta t$. The results are in agreement with \eqref{CM_method_error} where we observe global second-order accuracy for the test without submap decomposition and third-order accuracy with submap decomposition.\par

The first convergence test we perform consists of a Rossby-Haurwitz (RH) wave. These waves form exact time-dependent solutions to the Euler equations and play an important role in global atmospheric circulation \cite{pedlosky1987geophysical, haurwitz1940motion}. The RH wave is comprised of a stream function and vorticity each proportional to a spherical harmonic $Y_{\ell}^m$ rotating with constant phase speed $\nu = -2\Omega/\ell(\ell+1)$ \cite{haurwitz1940motion, neamtan1946motion}. The particular relative vorticity we use is given by 

\begin{equation}
\label{RHwave}
\zeta(\lambda, \theta,t) = 30 \cos(\theta)\sin^4(\theta) \cos(4(\lambda - \nu t)).
\end{equation}In addition to the vorticity \eqref{RHwave} we consider a non-rotating form of the RH wave in a rotated coordinate system where $(\lambda, \theta)$ are measured from an axes formed by applying a rotation by $\pi/3$ about the y-axis to the standard coordinate axes. This serves to demonstrate that the method suffers no constraints due to a choice of coordinate system since all computations are performed in Cartesian coordinates.\par 

The third convergence test consists of an initial relative vorticity distribution given by a Gaussian vortex of the form

\begin{equation}
\label{reversing_test}
\zeta_0(\bsym{x}) = 4\pi\text{\,exp}(-16 \norm{\bsym{x} - \bsym{x}_c}^2)\, ,
\end{equation}
restricted to the sphere, where the position of the centre of the vortex is taken to be $\bsym{x}_c = (1,0,0)$. We include only the conservation errors for this test since an analytic solution is not known. \par

We lastly consider an initial vorticity distribution given by a single zonal jet, forming an unstable steady solution to Euler's equations \cite{lorenz1972barotropic, tung1981barotropic}. We ran the test using a similar initial vorticity distribution as described in \cite{bosler2013particle}, given by 

\begin{equation}
\label{ZJ_vort}
\begin{aligned}
u(\lambda, \theta) &= \frac{\pi}{2}\text{exp}\left(-2\beta^2(1-\sin(\theta +\theta_c)) \right) \,,
\\
\zeta^{(zj)}_0(\lambda, \theta) &= \sin(\theta)(2\beta^2(\cos(\theta_c)\cos(\theta) - \sin(\theta_c)\sin(\theta)) + \cos(\theta))u(\lambda, \theta) \,,
\end{aligned}
\end{equation}
where $\theta_c$ is the centerline of the jet. The parameters are chosen to be $\beta = 12$ and $\theta_c = \pi/4$ with a final integration time of $T=0.5$. \par 

\begin{figure}[h]
\centering
\includegraphics[width = \linewidth, height = 10cm]{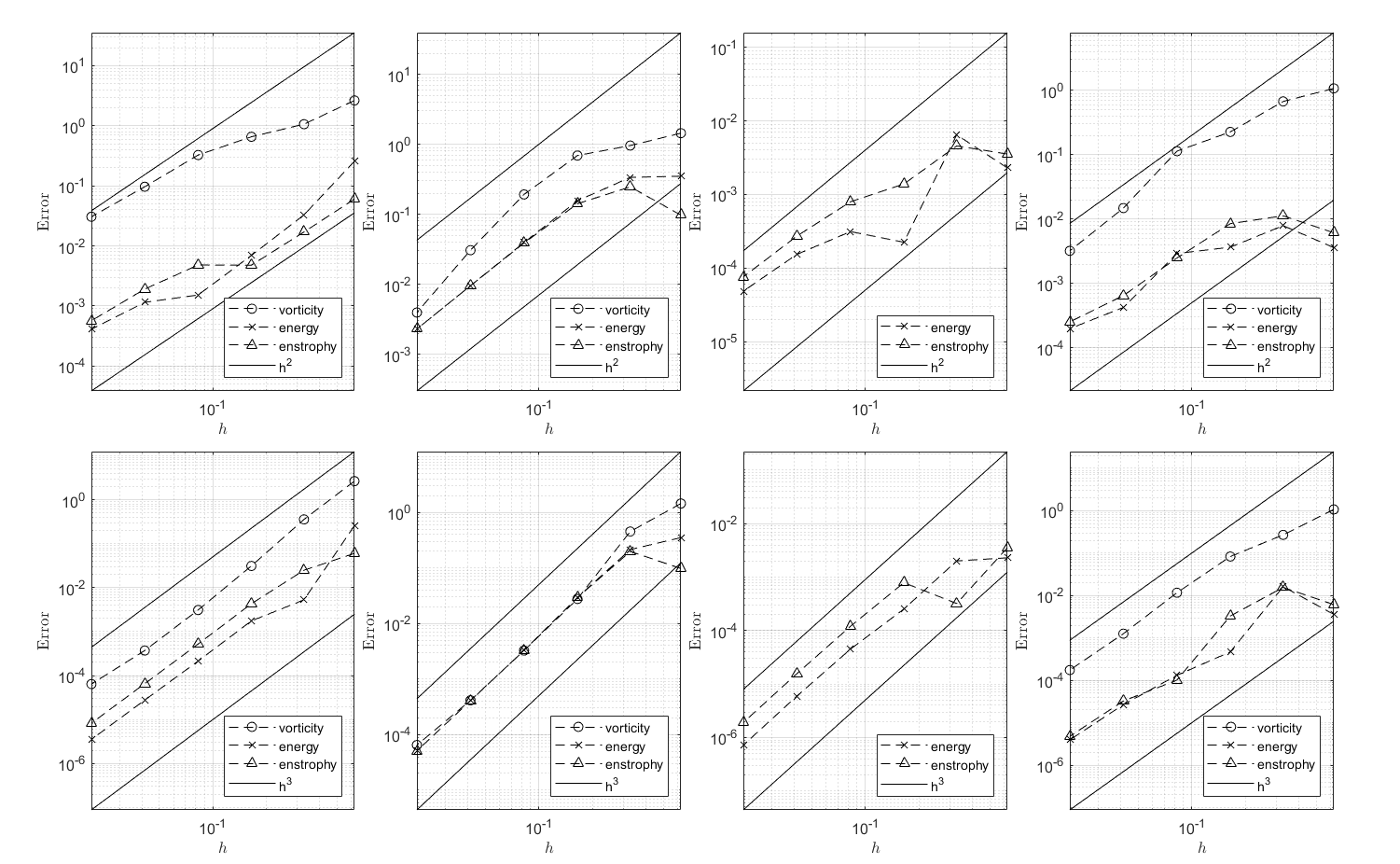}
\vspace{-0.5cm}
\caption{Convergence tests defined by the initial conditions \eqref{RHwave}, for the non-rotating \eqref{RHwave} in a rotated coordinate system (middle-left), Gaussian vortex \eqref{reversing_test}, and the zonal jet test case \eqref{ZJ_vort} from left to right. The top row consists of the tests performed without any submap decomposition and the tests in the bottom row used submap decomposition every $10$ time steps.}
\label{fig:euler_test_all}
\end{figure}

The results of the convergence tests are given in figure \ref{fig:euler_test_all}. We observe that in each test case without submap decomposition the convergence is globally second order accurate, and with submap decomposition we observe third order accuracy, affirming the theoretical predictions given in \eqref{error_estimates}. 

\section{Numerical Experiments}
\label{sec:numerical_experiments}
In this section we present numerical experiments designed to illustrate the resolution properties of the method. We perform simulations for initial vorticities which quickly transition into a turbulent flow with a multi-scale structure and measure their turbulent energy spectra at sub-grid scales. It is beyond the scope of this article to present an investigation through DNS with the proposed method for the late-time behaviour of the energy spectrum for incompressible, inviscid turbulence on the sphere. However, this problem has a noteworthy history and there remain a number of unresolved problems related to the organization of the non-linear evolution \cite{dritschel2015late,lindborg2022two}. In 1953, Fj{\o}rtoft recognized that the conservation of energy and enstrophy for two-dimensional incompressible turbulence on the sphere indicated a simultaneous cascade of enstrophy to small scales and energy to large scales \cite{fjortoft1953changes}. This theory was later brought into a quantitative form by Kraichnan, Leith, and Batchelor (KLB) \cite{kraichnan1967inertial,leith1968diffusion, batchelor1969computation}. The KLB theory of the dual cascade in the energy spectrum predicts a direct cascade proportional to $k^{-3}$ at small scales before the dissipative effects of viscosity and a $k^{-5/3}$ inverse cascade forming at the large scales \cite{boffetta2012two}. This prediction was made for the viscous case with an injection of energy, required to maintain a balance with dissipation. In the unforced and inviscid case, a mathematical description for the late-time behaviour of the energy spectrum is an open problem \cite{dritschel2015late}. An interesting recent work of Modin and Viviani \cite{modin2022canonical} provides insight into the mechanisms of a canonical scale separation for two-dimensional turbulence using the finite-mode approximation of Euler's equations on the sphere. \par   

The particular numerical experiments considered in this section consist of a randomly initialized vorticity distribution and a multiple perturbed zonal jet initial condition. We provide evidence of the ability to resolve scales in the vorticity that are beyond the computational grids defined by $L$ and observe the expected direct energy cascade at scales up to $\ell = 4096$. A heuristic explanation for the sub-grid resolution observed in these experiments is given in the \ref{AppendixA}. \par

\begin{figure}[h!]
\centering
\includegraphics[width = 11cm, height = 6cm]{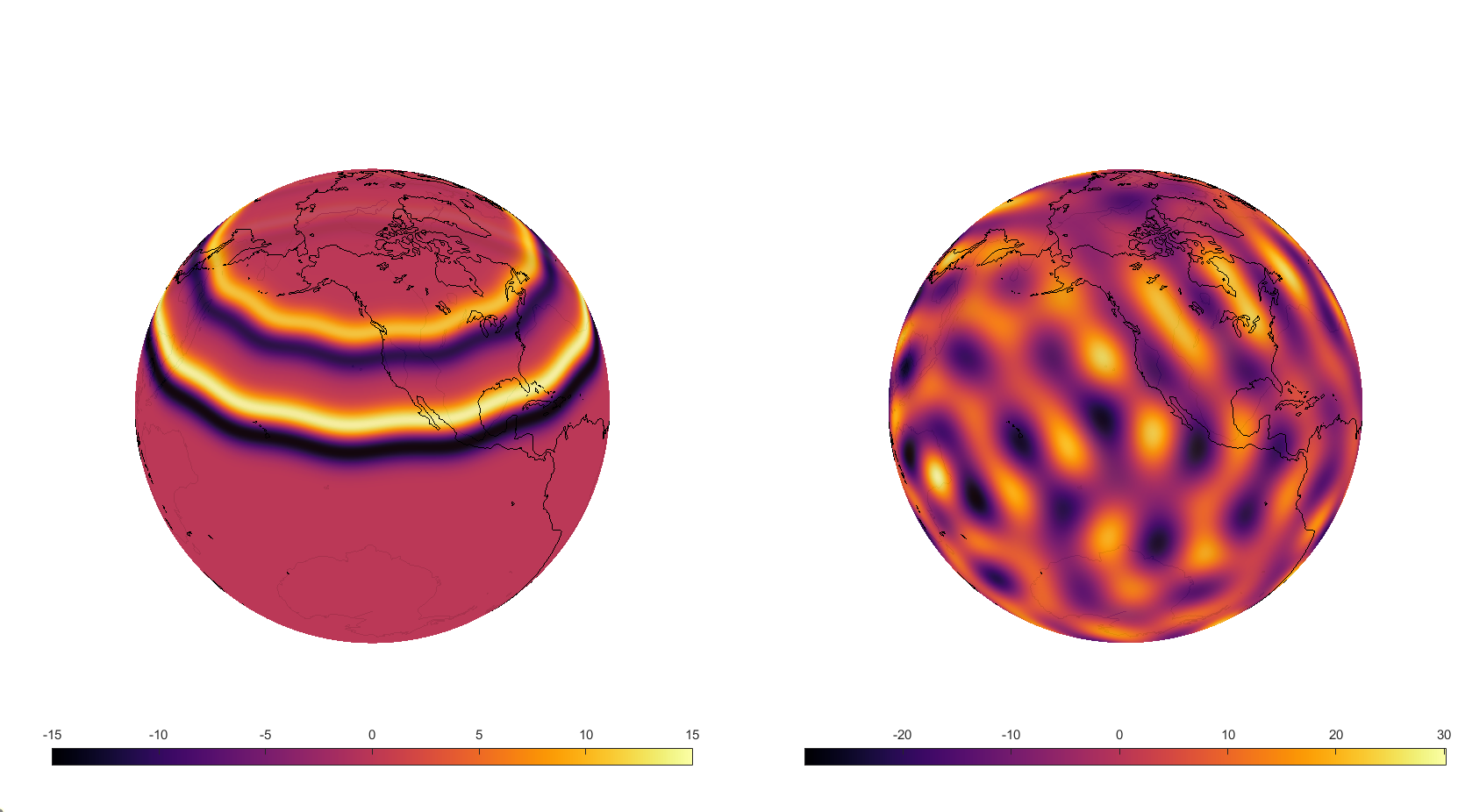}
\vspace{-0.5cm}
\caption{Initial Conditions for numerical experiments. Left: Perturbed multiple zonal jet initial condition. Right: Randomly initialized vorticity initial condition.}
\label{fig:ZJ_ics}
\end{figure}

\subsection{Multiple Perturbed Zonal Jets}

The modeling and simulation of the instability of zonal flows plays an important role in the understanding of stratospheric dynamics and the mechanisms behind sudden stratospheric warming events \cite{haynes2005stratospheric}. Numerical studies of these events necessitate the resolution of a large range of scales due to the production of large scale planetary waves along with the formation of vortex filaments, from which secondary vortices can form \cite{haynes2005stratospheric, bosler2013particle}. In an effort to demonstrate the capabilities of the method to simulate the complex vortex dynamics of a perturbed zonal jet, albeit in the idealized form of the rotating barotropic vorticity equations, we consider an initial condition formed by the sum of two zonal jets \eqref{ZJ_vort}. The perturbed centrelines of the jets are taken to be $\theta_{c1} = \pi/4 + 0.01 \cos(12\lambda)$ and $\theta_{c2} = 3\pi/8 + 0.01\cos(12 \lambda)$. As parameters of the simulation we used $T = 10$, $\Delta t = 1/1000$, the vorticity was evolved using a sampling defined by $L = 256$, and the $k=5$ level of refinement for the icosahedral discretization of the submaps. The simulation was performed with a submap decomposition every $10$ time steps.
\begin{figure}[H]
\includegraphics[width = \linewidth, height = 6cm]{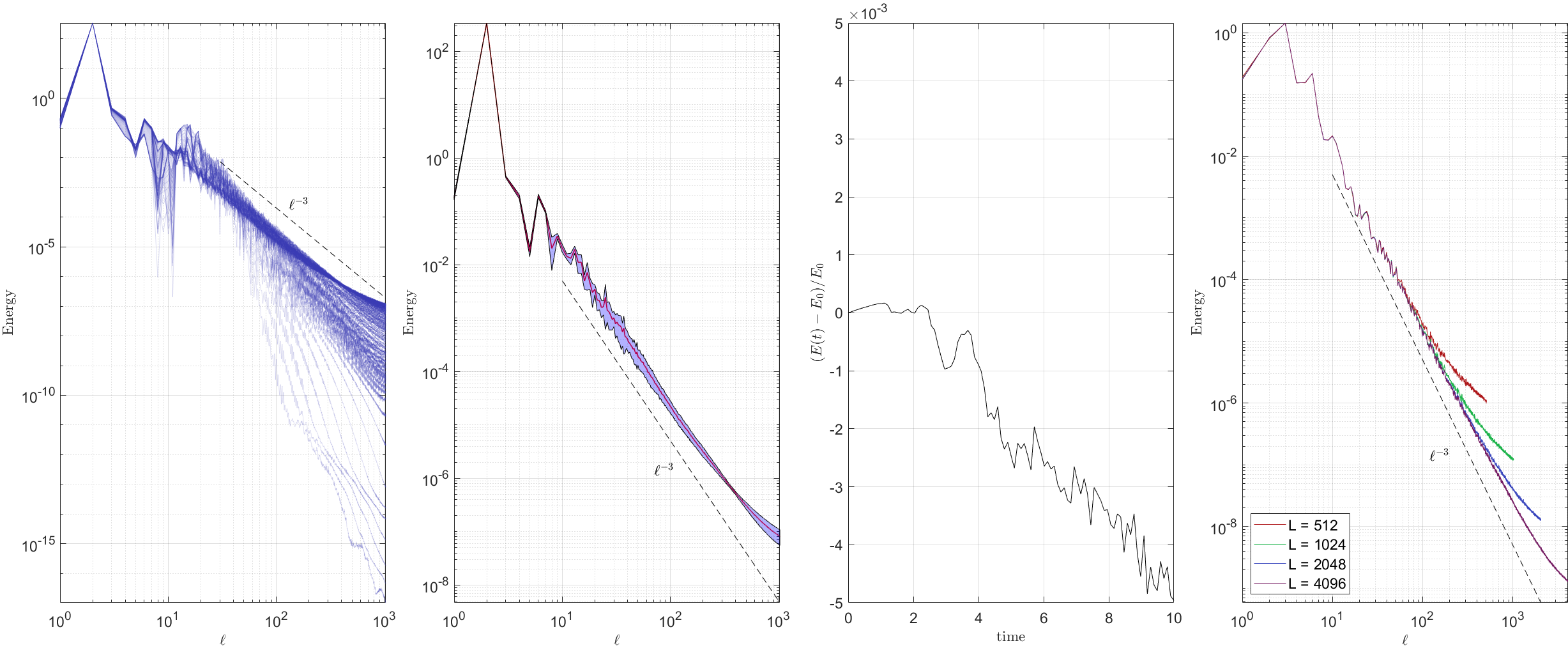}
\caption{Left: Evolution of the energy spectrum of the multiple zonal jet vorticity distribution with time increasing with the darkness of the lines. Middle-left: Mean and standard deviation of energy spectrum for $t = 5$ to $t = 10$. Middle-right: Change in energy in the first $\ell = 1000$ modes over time. Right: Energy spectrum of absolute vorticity at $t=10$ sampled up to band-limit $L = 4096$. Dashed black lines are proportional to $\ell^{-3}$.}
\label{fig:multi_jet_spectra}
\end{figure}

\begin{figure}
\centering
\includegraphics[width = 15cm, height = 11cm]{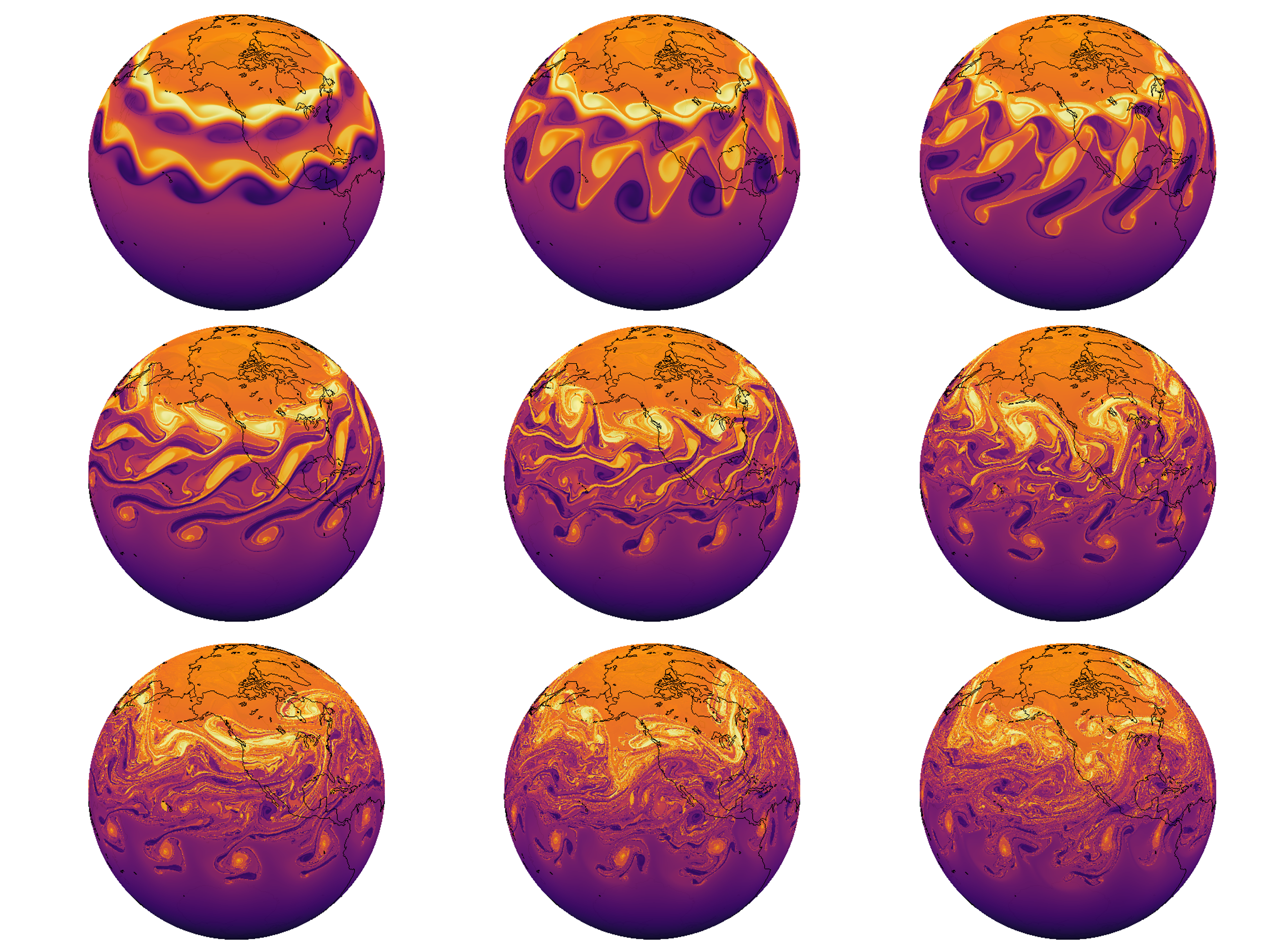}
\caption{Evolution of the rotating multiple zonal jet absolute vorticity from $t = T/10$ to $t = 9T/10$ in increments of $T/10$.}
\label{fig:multi_jet_evolution}
\end{figure}

The formation of distinct large scale vortices and connecting filaments is clearly observed after $t = T/2$. There is a subsequent strong mixing of all of these structures for the remainder of the simulation. Larger vortex structures begin to emerge and vortex filaments persist throughout the evolution resulting in a complex final vorticity distribution with a large range of spatial scales. We demonstrate the ability to retain vortex structures much finer than the computational grid in figure \ref{fig:multi_jet_evolution} with a zoom of the vorticity at the focal point $(\lambda, \theta) = (3.22055,1.1963)$  up to a window width of $2^{-12}$. The capacity of the method to capture sub-grid scale oscillations is observed at the final window width and could be refined further beyond what is depicted. \par

\begin{figure}
\centering
\includegraphics[width = \linewidth, height = 9.33cm]{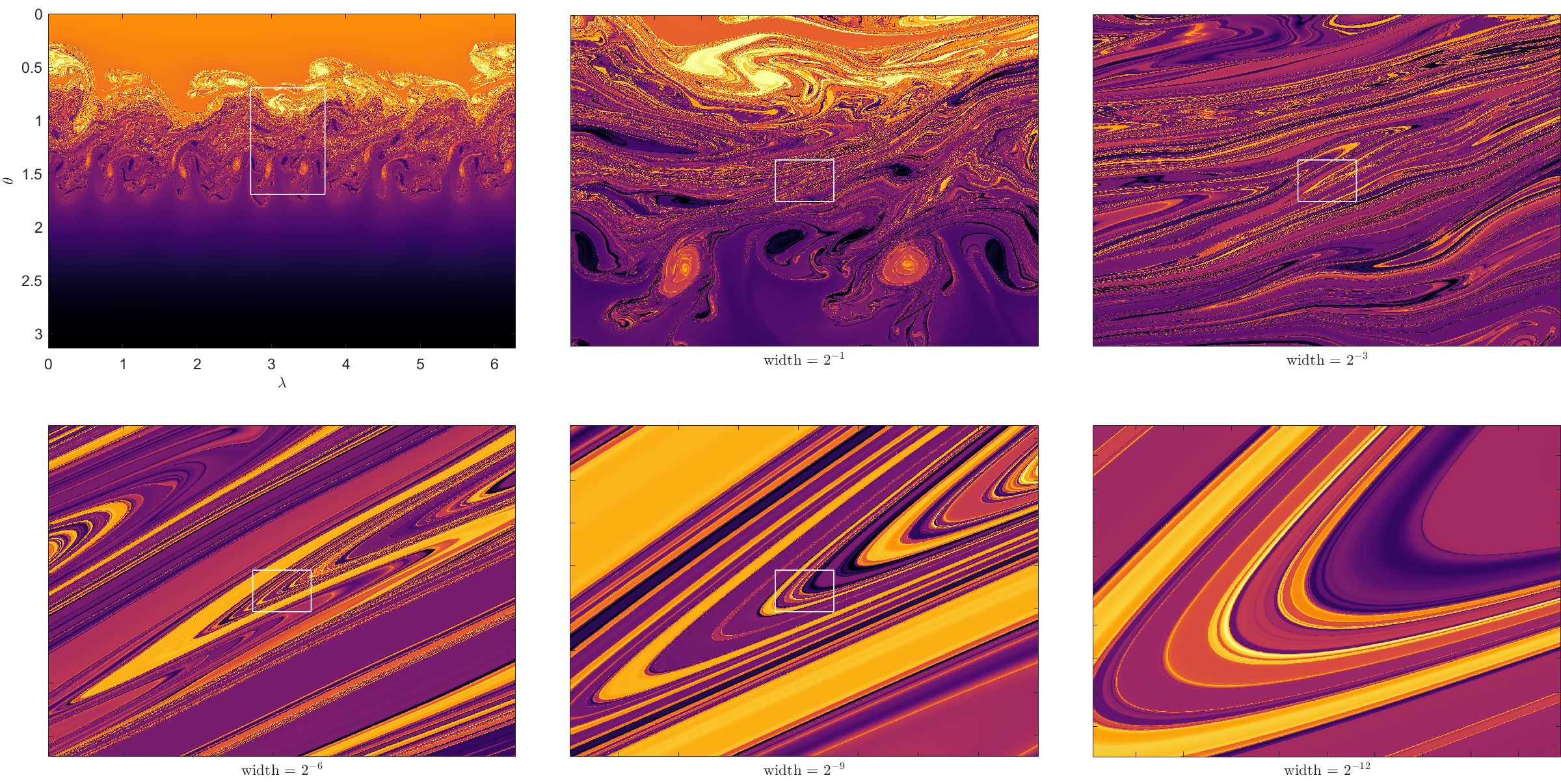}
\vspace{-0.5cm}
\caption{Zoom into vortex structures up to a width of $2^{-12}$ of \ref{fig:multi_jet_evolution} at $t=10$.}
\end{figure}

The energy spectrum and conservation properties of the simulation are shown in figure \ref{fig:multi_jet_spectra}. In the middle-right panel we have plotted the normalized error in the energy in the first $ \ell \leq 1000$ frequencies over time and observe the energy to decrease indicating a downscale transfer of energy. The cascade of energy towards small scales can be seen over the course of the evolution and a distinct $\ell^{-3}$ scaling forms in the energy spectrum over the last half of the simulation. This scaling is shown to carry down towards wave numbers of $\ell = 4096$.  The inverse cascade scaling of $\ell^{-5/3}$ is however not observed and instead the $\ell^{-3}$ scaling becomes more prominent at lower wavelengths.

\subsection{Randomly Initialized Vorticity}

\begin{figure}[h]
\includegraphics[width = \linewidth, height =5.5cm]{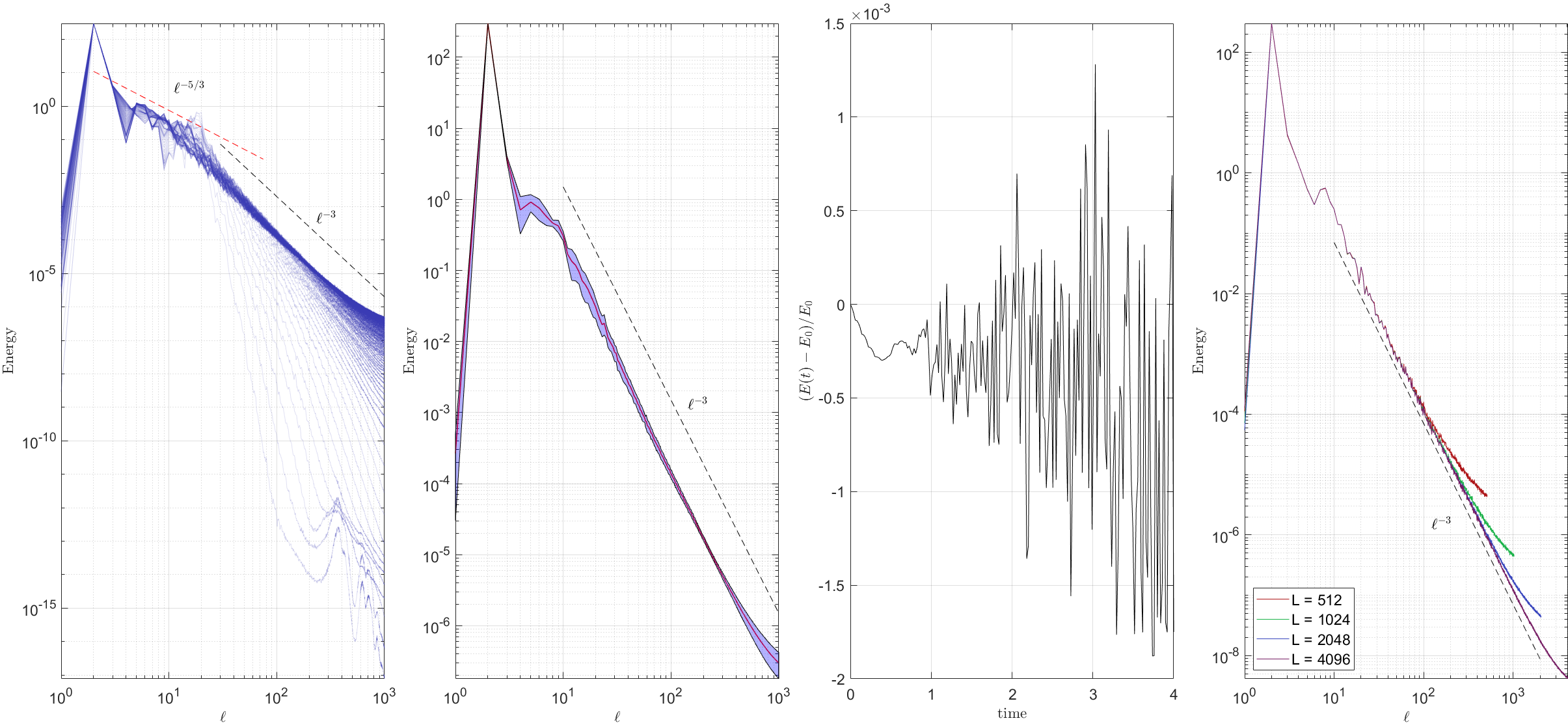}
\caption{Left: Evolution of the energy spectrum of the randomly initialized vorticity with time increasing with the darkness of the lines. Middle-left: Mean and standard deviation of energy spectrum for $t = 400\Delta t$ to $t = T$. Middle-right: Change in energy in the first $\ell = 1000$ modes over time. Right: Energy spectrum of absolute vorticity at $t=4$ sampled up to band-limit $L = 4096$. Dashed black lines are proportional to $\ell^{-3}$ and the dashed red line is proportional to $\ell^{-5/3}$.}
\label{fig:random_vorticity_spectra}
\end{figure}

We performed a simulation for a vorticity distribution initialized as a sum of the first $\ell = 20$ with all $m$ modes filled such that the vorticity was real and with amplitudes sampled from a uniform distribution over the interval $[-5,5]$.  As parameters of the simulation we used  $T= 4$, $\Delta t = 1/1000$, a band-limit $L = 256$ for the velocity field and the $k = 6$ level of refinement for the icosahedral discretization of the submaps. The simulation is performed with a fixed remapping strategy over every $20$ time steps. We observe the ability to retain scales beyond the coarser computational grid of the velocity field by computing the energy spectrum up to a band-limit of $L = 1000$ in figure \ref{fig:random_vorticity_spectra}. The expected energy cascade of $E(\ell,t) \sim \ell^{-3}$ is observed as $t \to T$ and the inverse cascade scaling of $\ell^{-5/3}$ emerges in intermediate times. As time progresses the range where this inverse cascade formed shrinks and is replaced by a scaling closer to $\ell^{-3}$. The energy behaves in a more oscillatory manner over time than for the multiple zonal jet simulation, indicating both a downscale and upscale transfer of energy in the first $1000$ modes.

\begin{figure}[H]
\centering
\includegraphics[width = 15cm, height = 6.5cm, trim=0cm 4cm 0cm 2.5cm, clip]{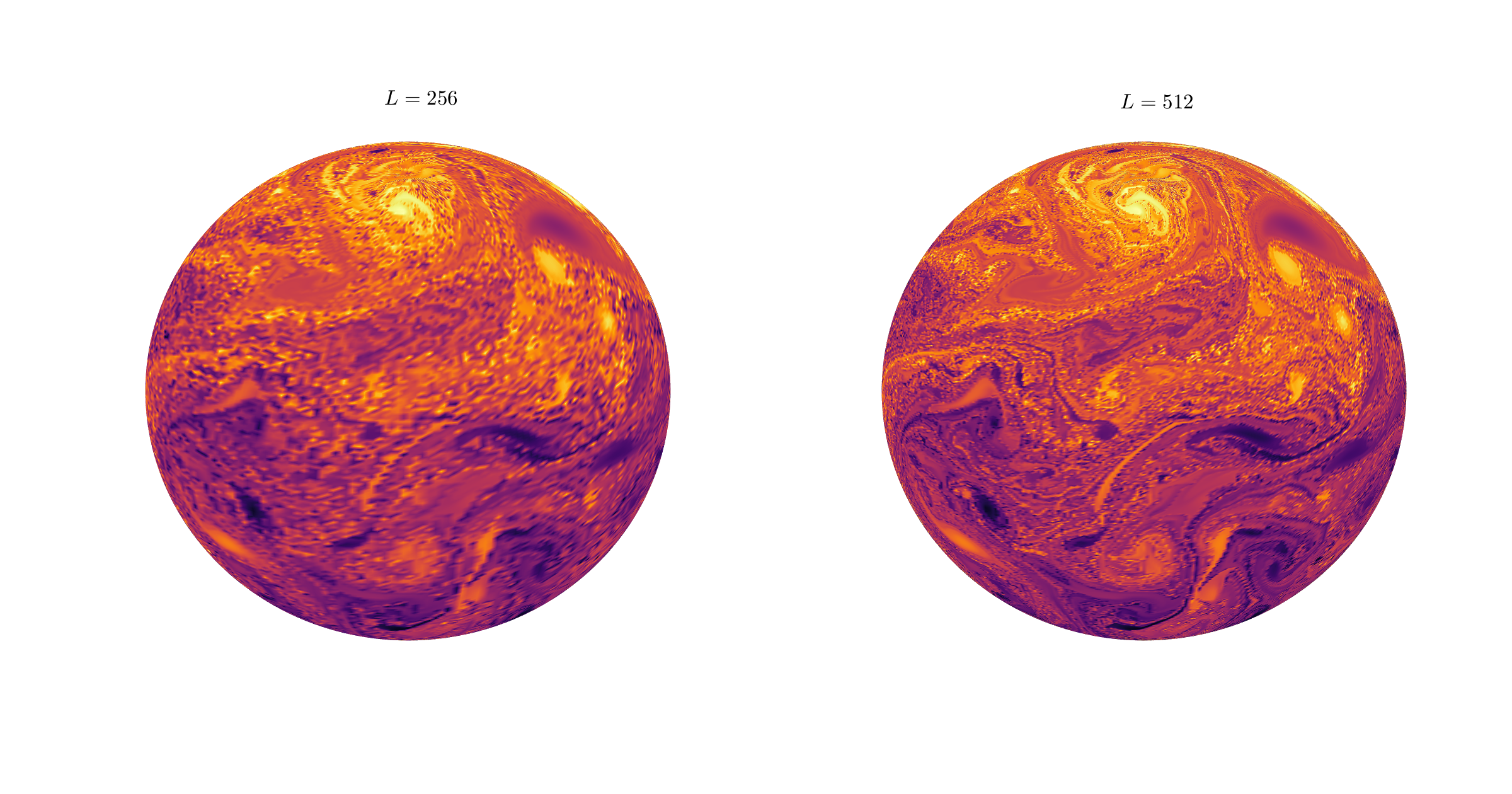}
\includegraphics[width = 15cm, height = 6.5cm, trim=0cm 4cm 0cm 2.5cm, clip]{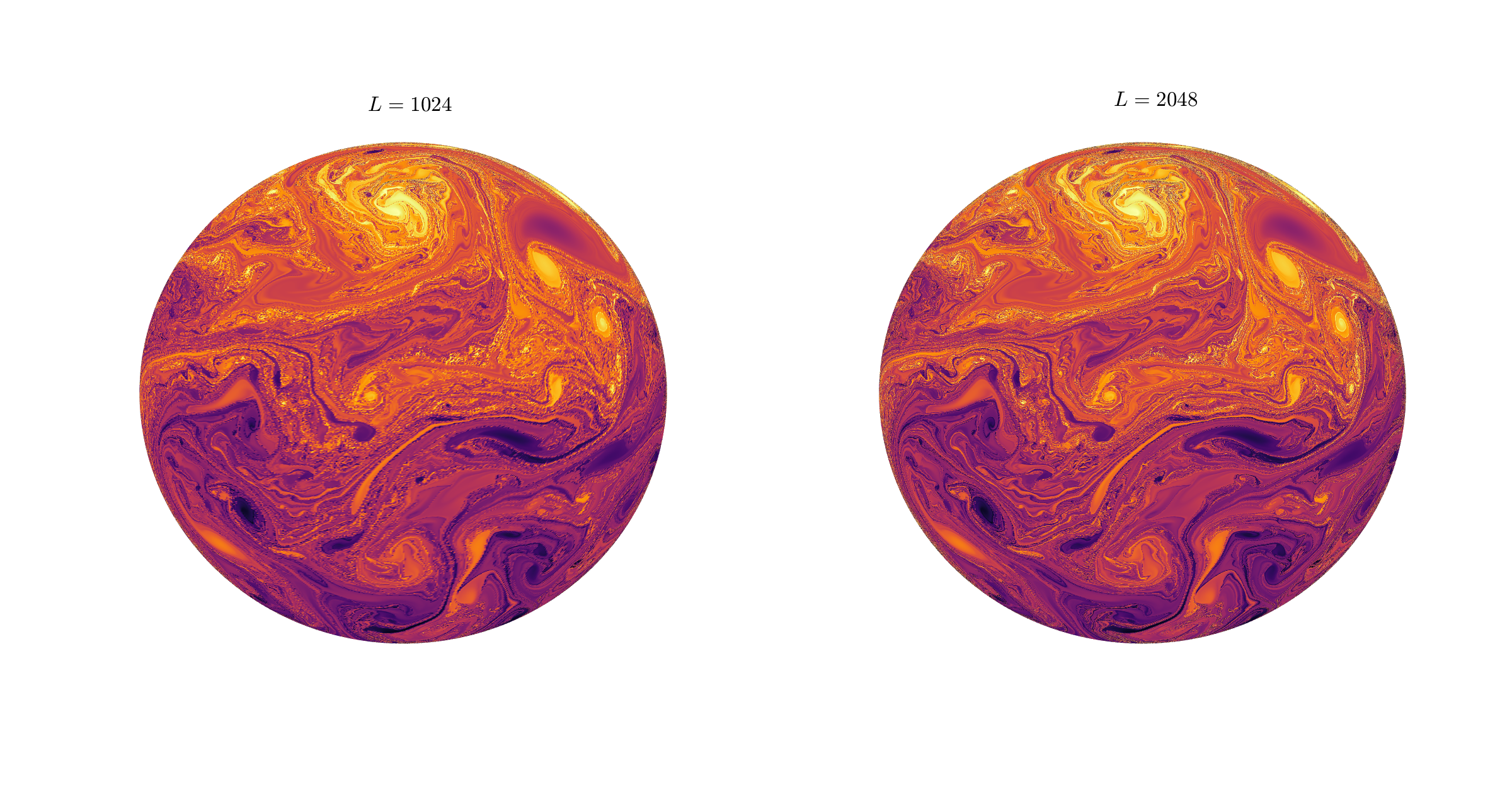}

\caption{Upsampling of the solution at $t=T$ for the randomly initialized vorticity distribution test case. A sampling on the grid points \eqref{grid_points} was performed for increasing band-limit $L$ and a distinct sharpening of the solution is observed with the anticipated energy cascade as indicated of the rightmost panel of figure \ref{fig:random_vorticity_spectra}.}
\end{figure}

\section{Conclusion and Outlook}

A semi-Lagrangian characteristic mapping method for incompressible hydrodynamics on a rotating sphere was presented. The method utilizes a spatio-temporal discretization of the inverse map generated by the Eulerian velocity field as a composition of sub-interval flows and the vorticity is evolved through pullback of the initial condition with this map. Each submap is discretized using an embedding-based approach for manifold-valued data approximation with piecewise spherical spline interpolation, extending the techniques developed in \cite{taylor2023projection} to a non-linear advective setting. The method was described in detail and error estimates were provided and validated using a number of standard test cases in section \ref{numerical_verification}, indicating global third-order accuracy in the supremum norm. Numerical experiments illustrating the unique resolution properties gained through the spatio-temporal discretization of the inverse map were performed in section \ref{sec:numerical_experiments}. The ability to reproduce the expected turbulent energy cascades and resolve vortex structures at sub-grid scales was demonstrated. \par 

We have differed a number of lines of investigation which can build on the techniques outlined here to future work. We are seeking to better understand and provide more rigorous mathematical justification for the resolution properties of the method. An analysis of the separation of the scales driving the evolution, filtered by the sampling with the vorticity, and the scales represented through the composition could elucidate an optimal range for the parameters of the method such as the spatial discretization of the map, the filter parameter $L$, and the number of compositions. Furthermore, the presented discretization of the inverse map can be applied to a broader class of equations on the diffeomorphism group and extended to other manifolds using modifications of the embedding-based approach. Additionally, since we have formulated the method on arbitrary triangulations it is directly amenable to the techniques of $h$- and $r$-adaptivity. Investigating the approximation capabilities of the submap decomposition using adaptive mesh refinement techniques is warranted. \par 

Beyond these numerical investigations, developing the CM method techniques to incorporate the effects of compressiblity for the shallow-water equations, along with other physical effects of advected parameters by the fluid are the subject of our current research. Since the method relies upon the transport structure of Euler's equations, it is not directly applicable to flows with diffusive processes such as the Navier-Stokes equations. Incorporating the effects of viscosity along with sub-grid scale thermodynamic processes will be important extensions of the method. Given the capability of representing a large range of spatial scales through the use of the spatio-temporal discretization \eqref{remapping}, we believe that computational and theoretical advancements of the method presented here will yield useful tools for the simulation and study of geophysical fluid dynamics.    



\section*{Acknowledgments}

The work of S.T. was partially supported by the NSERC CGS-D program.  The work of J-C.N. was partially supported by the NSERC Discovery Grant program and the Agence Nationale de la Recherche (ANR), grant ANR-20-CE46-0010-01. The authors would like to thank Xi-Yuan Yin of {\'E}cole Centrale de Lyon and S.T. would like to thank Tim Whittaker of UQAM for helpful discussions.

\appendix
\section{Lie advection of vorticity}
\label{remark1}
The Lie advection of the vorticity can be observed using the language of differential forms. Denote $(\cdot)^{\flat}: \mathfrak{X}(\mathbb{S}^2) \to \Omega^1(\mathbb{S}^2)$ as the flat operator, mapping vector fields into differential one-forms via the metric and let $\star: \Omega^k(\mathbb{S}^2) \to \Omega^{2-k}(\mathbb{S}^2)$ be the Hodge star operator. The planetary vorticity can be defined with respect to a potential $\bsym{R}(x) \in \mathbb{R}^3$ as the two-form $d\bsym{R}^{\flat} = 2\star \bsym{\Omega}^{\flat} \in \Omega^2(\mathbb{R}^3)$ restricted to the sphere. The covariant form of the incompressible Euler equations on a rotating sphere are 
 \begin{equation}
 \label{Euler_contra}
 \partial_t \bsym{u}^{\flat} + \mathcal{L}_{\bsym{u}} \bsym{u}^{\flat}  + \bsym{i}_{\bsym{u}}(d\bsym{R}^{\flat}) = - d(p-|\bsym{u}|^2/2 )\, , \quad \text{div}(\bsym{u})\mu = \mathcal{L}_{\bsym{u}}\mu = 0,
 \end{equation}where $p: \mathbb{S}^2 \to \mathbb{R}$ is the pressure and $i_{\bu}$ is the interior product with $\bu$, and $\mu$ is the Riemannian volume form. The absolute vorticity is given by $\bsym{\omega} = \bsym{\zeta} + d\bsym{R}^{\flat}$ where $\bsym{\zeta} = d\bu^{\flat}$ is the relative vorticity. Note then that since $M = \mathbb{S}^2$ for every $\bsym{\sigma} \in \Omega^2(\mathbb{S}^2)$ there is an associated $\sigma \in \Omega^0(\mathbb{S}^2)$ such that $\star \sigma = \bsym{\sigma}$. Expressing $\bsym{\omega} = \star \omega = \omega \mu$, taking the exterior derivative of \eqref{Euler_contra}, using the incompressibility constraint, along with Cartan's formula, we see that
\begin{equation}
\label{pullback_vort}
 (\partial_t  + \mathcal{L}_{\bsym{u}})\bsym{\omega}  = (\partial_t\omega + \mathcal{L}_{\bu}\omega)\mu = 0 \implies \omega(t) = \omega_0 \circ \varphi_{[t,0]}\,.
\end{equation}
The velocity field is then recovered from the stream function $\bsym{\psi} \in \Omega^2(\mathbb{S}^2)$ via $\bsym{u}^{\flat} = \delta\bsym{\psi}$ where $\delta = -\star d\star$ is the codifferential operator on $\Omega^2(\mathbb{S}^2)$. In vector form we have that $\bsym{u} = (\delta\bsym{\psi})^{\sharp} = -(\star d \star\star \psi)^{\sharp} =  (\star d\psi)^{\sharp} =  -\nabla^{\perp}\psi$. The relative scalar field vorticity $\zeta$ and scalar stream function $\psi$ are then related by the Poisson equation 
\begin{equation}
\label{Poisson_stream}
\zeta = \star^{-1}\star\zeta = -\star\bsym{\zeta} = -\star(d\delta \bsym{\psi}) = \star d(\star d\star\star\psi) = \Delta \psi \, ,
\end{equation}where $\Delta = -\star d(\star d) $ is the Hodge Laplacian on scalar functions. Note that the Hodge Laplacian defining the Poisson equation \eqref{Poisson_stream} differs from the surface Laplacian by a factor of $-1$: defining $\Delta_{g} = \text{div}(\nabla_{g})$ with respect to the metric and using $\text{div}(\bsym{u}) = -\delta\bsym{u}^{\flat}$, we see that $\Delta_{g}f = \text{div}(df)^{\sharp} = -\delta d f = \star d \star d f$ for a scalar function $f$. The solution of \eqref{Poisson_stream} is determined up to a constant and uniqueness is recovered by imposing a zero mean condition on the stream function. 

\section{Spherical harmonics and the angular momentum operator}
\label{sph_harmonic_appendix}

The space of spherical harmonics of degree $\ell$, denoted  $\mathcal{Y}_\ell$, is formed by the restriction of the space of harmonic homogeneous trivariate polynomials, that is 
\begin{equation}
\mathcal{Y}_{\ell} = \left\{ \left.p\right\vert_{\mathbb{S}^2} \,:\,\, p \in \mathcal{H}_{\ell}\,,\,  \Delta p = 0 \right\} \,,
\end{equation}where $\Delta$ is the Laplacian in Euclidean space. In spherical coordinates $(\lambda, \theta) \in [0, 2\pi) \times [0,\pi]$, we can write the basis functions for $\mathcal{Y}_{\ell}$, normalized in $L^2(\mathbb{S}^2)$, as
\begin{equation}
\label{sph_harms}
    Y_{\ell}^m(\lambda, \theta) = (-1)^m\sqrt{\frac{2\ell + 1}{4\pi}\frac{(\ell - m)!}{(\ell + m)!}} P_{\ell}^m(\cos(\theta))e^{im\lambda} \, ,
\end{equation}
where $|m| \leq \ell$ and $P^m_{\ell}(x)$ are the associated Legendre functions.  The orbital angular momentum operator is used to compute the rotated spherical gradient relating the stream function to the velocity field. This operator acts on the spherical harmonic basis as
\begin{equation}
\label{ang_op}
\bsym{L}_{\hat{\bsym{n}}}Y_{\ell,m}(\bsym{x}) = -i \hat{\bsym{n}} \times \nabla Y_{\ell,m}(\bsym{x}) \,.
\end{equation}
The operation \eqref{ang_op} can be written as block diagonal matrix on the coefficients.  We express $\bsym{L}$ as a linear combination of the infinitesimal rotations about the standard Cartesian basis vectors. The components of \eqref{ang_op} in the $x$ and $y$ directions are given by

\begin{equation}
\label{pm_ops}
\bsym{L}_x = \frac{1}{2}\left(\bsym{L}_+ - \bsym{L}_{-}\right) \, , \quad \bsym{L}_y = \frac{i}{2}\left(\bsym{L}_{-} - \bsym{L}_{+}\right) \, .
\end{equation}where $\bsym{L}_{\pm}$ are the raising and lowering operators. These operators and the $z$-component of the angular momentum act on the spherical harmonic basis functions as

\begin{equation}
\label{L_operator}
\bsym{L}_{\pm} Y_{\ell}^m = \sqrt{\ell(\ell + 1) - m(m\pm 1)}Y^{m\pm 1}_\ell \, , \quad \bsym{L}_z Y_{\ell}^m = mY_{\ell}^m \, .
\end{equation}

\section{Convergence of the method}
\label{convergence_proof}
\begin{pf}{(of Theorem \ref{theorem1})}
We include the remapping estimate \eqref{remapping_error} with an analogous proof as given in \cite{yin2021characteristic}. Introducing the map $\tilde{\varphi}_{[T,0]}$ resulting from the modified equation and splitting the error, it suffices to estimate the error between $\varphi_{[T,0]}$ and $\tilde{\varphi}_{[T,0]}$. Considering the error to $t_n$, using the Lipschitz continuity of $\tilde{\varphi}$ we have that  

\begin{equation}
\begin{aligned}
\norm{\varphi_{[t_n,0]}- \tilde{\varphi}_{[t_n,0]}}_{C^{0,\alpha}} &\leq \norm{(\varphi_{[t_{n-1},0]}- \tilde{\varphi}_{[t_{n-1},0]}) \circ \varphi_{[t_n, t_{n-1}]}}_{C^{0,\alpha}} 
\\
&+ \norm{\tilde{\varphi}_{[t_{n-1},0]}\circ \varphi_{[t_n,t_{n-1}]} - \tilde{\varphi}_{[t_{n-1},0]}\circ\tilde{\varphi}_{[t_n,t_{n-1}]}}_{C^{0,\alpha}} 
\\
& \leq \norm{\varphi_{[t_{n-1},0]} - \tilde{\varphi}_{[t_{n-1},0]}}_{C^{0,\alpha}} + C_{n-1}\norm{\varphi_{[t_n,t_{n-1}]} - \tilde{\varphi}_{[t_n,t_{n-1}]}}_{C^{0,\alpha}} \,.
\end{aligned}
\end{equation}
where $C_{n-1}$ is the Lipschitz constant of $\tilde{\varphi}_{[t_{n-1},0]}$. The second term can be estimated as

\begin{equation}
\begin{aligned}
\norm{\varphi_{[t_n,t_{n-1}]} - \tilde{\varphi}_{[t_n,t_{n-1}]}}_{C^{0,\alpha}} &\leq \int_{t_{n-1}}^{t_n} \norm{\bsym{u}(\varphi_{[s,t_{n-1}]},s) - \tilde{\bsym{u}}(\tilde{\varphi}_{[s,t_{n-1}]},s)}_{C^{0,\alpha}} ds
\\
& \leq \int_{t_{n-1}}^{t_n}\norm{(\bsym{u}-\tilde{\bsym{u}})(\varphi_{[s,t_{n-1}]},s)}_{C^{0,\alpha}} + K(s)\norm{\varphi_{[s,t_{n-1}]} - \tilde{\varphi}_{[s,t_{n-1}]}}_{C^{0,\alpha}} ds \,.
\end{aligned}
\end{equation}
where $K(s)$ is the Lipschitz constant of $\tilde{\bu}(s)$. Then by applying Gr{\"o}nwall's lemma we can say that  

\begin{equation}
\norm{\varphi_{[t_n,t_{n-1}]} - \tilde{\varphi}_{[t_n,t_{n-1}]}}_{C^{0,\alpha}} = \mathcal{O}(\Delta t \sup_{t \in [t_{n-1},t_n]} \norm{\bsym{u}(t)-\tilde{\bsym{u}}(t)}_{C^{0,\alpha}}) \,.
\end{equation}
Over the interval $[t_{n-1},t_n]$ the velocity field is extrapolated using the $p$ previous iterations, allowing us to write

\begin{equation}
\sup_{t \in [t_{n-1},t_n]}\norm{\bsym{u}(t)-\tilde{\bsym{u}}(t)}_{C^{0,\alpha}} = \mathcal{O}(\Delta t^{p} + \norm{\bsym{u}(t_{n-1})-\tilde{\bsym{u}}(t_{n-1})}_{C^{0,\alpha}}) \,.
\end{equation}
We can bound the error in the velocity field by the error in the map using a Schauder estimate for the Laplacian \cite{gilbarg1977elliptic} of the form

\begin{equation}
\label{velocity_bound}
\begin{aligned}
\norm{\bsym{u}(t_{n-1})-\tilde{\bsym{u}}(t_{n-1})}_{C^{0,\alpha}} &\lesssim \norm{\psi(t_{n-1}) - \tilde{\psi}(t_{n-1})}_{C^{2,\alpha}} 
\\
&\lesssim \norm{\omega(t_{n-1}) - \tilde{\omega}(t_{n-1})}_{C^{0,\alpha}}
\lesssim \norm{\varphi_{[t_{n-1},0]} - \mathcal{X}_{[t_{n-1},0]}}_{C^{0,\alpha}} \,.
\end{aligned}
\end{equation}Combining these estimates and splitting the error in \eqref{velocity_bound} again with $\tilde{\varphi}_{[t_{n-1},0]}$ we get that

\begin{equation}
\begin{aligned}
\norm{\varphi_{[t_n,0]}- \tilde{\varphi}_{[t_n,0]}}_{C^{0,\alpha}} &\lesssim \norm{\varphi_{[t_{n-1},0]}- \tilde{\varphi}_{[t_{n-1},0]}}_{C^{0,\alpha}} + \mathcal{O}(\Delta t (\Delta \tau \min(h^3\Delta t^{-1}, h^2) +  \Delta t^s) + \Delta t^{p+1})
\end{aligned}
\end{equation}
Setting $t_n = T$ and iterating this argument for the previous maps implies the desired result \eqref{CM_method_error}. 
\qed
\end{pf}

\section{The resolution of the submap decomposition}
\label{AppendixA}
The finest scale contained in the approximation of the vorticity is dictated by the spatial discretization of the submaps and the number of composition forming the decomposition \eqref{remapping}. Since this approximation can support a global representation of fine scales, we can associate a notion of resolution to an effective band-limit of the discretization in terms of spherical harmonics. In this appendix we indicate how the discretization is able to represent exponentially fine scales  through simple rules for compositions of global polynomials. \par 

If we consider initial vorticities which are band-limited then the increase in support in frequency space over time is dictated by the scales generated in the inverse map since
\begin{equation}
\label{effective_resolution}
\omega_0 \circ \mathcal{X}_{[t,0]}= \sum_{\ell = 0}^{L_0}\sum_{|m|\leq \ell} \hat{\omega}^0_{\ell,m}Y_{\ell}^m \circ \mathcal{X}_{[t,0]}  =  \sum_{\ell' = 0}^{L_e(t)}\sum_{|m'|\leq \ell'} \sum_{\ell = 0}^{L_0}\sum_{|m| \leq \ell} \hat{\omega}^0_{\ell,m} b_{\ell',m'}^{\ell,m}(t) Y_{\ell'}^{m'}  + \mathcal{R}(\mathcal{X}_{[t,0]})\,,
\end{equation}
where the coefficients are given by $b_{\ell',m'}^{\ell,m}(t) = \inprod{Y_{\ell}^m \circ \mathcal{X}_{[t,0]}}{Y_{\ell'}^{m'}}$ and $\mathcal{R}(\mathcal{X}_{[t,0]})$ is a remainder term depending on the spatial discretization of the submaps. We call $L_e(t)$ the effective band-limit since we do not necessarily have access to all of these scales globally, nor is the vorticity ever expanded up to this wavenumber. Rather, the spatial representation of the vorticity admits an effective representation in terms of spherical harmonics which we can use to quantify the frequency content present in this approximation. Since the method gives a spatial, rather than spectral, approximation of the vorticity, we believe a more apt description would be given using a multi-resolution analysis incorporating both space and scale. \par
In order to simplify the discussion, we consider a discretization of the map using spherical harmonics which eliminates the remainder term. The notion of effective band-limit then becomes clear with an elementary analysis of composition with projected spherical polynomials. Let the space of projected spherical vector polynomials be
\begin{equation}
\mathcal{P}(\mathcal{B}_d^{3}) = \left\{\mathcal{P}(\bsym{p}) \,:\, \bsym{p} \in \mathcal{B}_d^3\,,\,\,  \bsym{p}(v) \in \mathbb{R}^3\setminus\{0\} \,\,\,  \forall v \in \mathbb{S}^2\right\} \,.
\end{equation} 
The homogeneous polynomial maps $F \in \mathcal{P}(\mathcal{B}_d^3)$ admit degree $d$ extensions $\bar{F} \in \mathcal{H}_d^3$ such that 

\begin{equation}
\bar{F}(\bsym{x}) = \norm{\bsym{x}}^d F(\bsym{x}/\norm{\bsym{x}}) \,, 
\end{equation}
for all $\bsym{x} \in \mathbb{R}^3\setminus\{0\}$. Using the homogeneous extensions of two maps $F_1, F_2 \in \mathcal{P}(\mathcal{B}_d^3)$ we see that

\begin{equation}
\mathcal{P}(F_2) \circ \mathcal{P}(F_1)(\bsym{x}) = \frac{\bar{F}_2(\bar{F}_1(\bsym{x})/\norm{\bar{F}_1(\bsym{x})})}{\norm{\bar{F}_2(\bar{F}_1(\bsym{x})/\norm{\bar{F}_1(\bsym{x})})}} = \frac{\norm{\bar{F}_1(\bsym{x})}^{-d}\bar{F}_2(\bar{F}_1(\bsym{x}))}{\norm{\bar{F}_1(\bsym{x})}^{-d}\norm{\bar{F}_2(F_1(\bsym{x}))}}  = \mathcal{P}(\bar{F}_2 \circ \bar{F}_1)(\bsym{x}) \,.
\end{equation}Noting that $\mathcal{H}_{\ell} \circ (\mathcal{H}_d^3) \subset \mathcal{H}_{d\cdot \ell} \Rightarrow \mathcal{H}_{l} \circ (\mathcal{B}_d^3)\subset \mathcal{B}_{d\cdot \ell}$ we see 
\begin{equation}
\label{property1}
\mathcal{P}(\mathcal{B}_d^3)\circ\mathcal{P}(\mathcal{B}_d^3) \subset \mathcal{P}((\mathcal{H}_d \circ \mathcal{B}_d^3)^3) \subset \mathcal{P}(\mathcal{B}_{d^2}^3) \,.
\end{equation} Using this inclusion we can deduce properties of the composition of projected vector spherical polynomials from the properties of $\mathcal{B}_d^3$. 

\begin{prop}
\begin{equation}
\mathcal{Y}_{\ell} \circ \mathcal{P}(\mathcal{B}_d^3) \subset \mathcal{B}_{d\cdot \ell} \,.
\end{equation}
\end{prop}
\begin{pf}
Let $\bar{Y}_{\ell}^m \in \mathcal{H}_{\ell}$ be the homogeneous extension of degree $\ell$ of the spherical harmonic $Y_{\ell}^m$ and $\bar{\mathcal{X}}$ the homogeneous extension of degree $d$ of $\Xhi \in \mathcal{P}(\mathcal{B}_d^3)$.  We have that for all $\bsym{x} \in \mathbb{R}^3\setminus\{0\}$ and $\lambda > 0$   
\begin{equation*}
\bar{Y}_{\ell}^m\circ \bar{\Xhi}(\lambda\bsym{x}) = \bar{Y}_{\ell}^m \left(\lambda^d\norm{\bsym{x}}^d\mathcal{X}(\bsym{x}/\norm{\bsym{x}})\right) = \lambda^{\ell \cdot d}\norm{\bsym{x}}^{\ell \cdot d}Y_{\ell}^m \circ \Xhi(\bsym{x}/\norm{\bsym{x}}) \,,
\end{equation*}and therefore $\bar{Y}_{\ell}^m \circ \bar{\Xhi} \in \mathcal{H}_{\ell \cdot d}$ from which the claim follows by restriction to the sphere. \qed
\end{pf}
In \cite{alfeld1996bernstein} it was observed that $\mathcal{B}_{\ell}$ admits the following decomposition in terms of the spherical harmonic spaces

\begin{equation}
\label{Bd_decompostion}
\mathcal{B}_{\ell} = \begin{cases}
\mathcal{Y}_0\oplus \mathcal{Y}_{2} \oplus \dots \oplus \mathcal{Y}_{2d} \,, &\text{ if } 2d = \ell
\\
\mathcal{Y}_1\oplus \mathcal{Y}_{3} \oplus \dots \oplus \mathcal{Y}_{2d-1} \,, &\text{ if } 2d-1 = \ell \,.
\end{cases}
\end{equation}
As a consequence, the composition of a degree $\ell$ spherical harmonic with a projected spherical polynomial of degree $d$ will admit an expansion in spherical harmonics of only even or odd degree based on the parity of $d \cdot \ell$. General spherical band-limited functions of degree $\ell$ are in the direct sum $\mathcal{B}_{\ell} \oplus \mathcal{B}_{\ell-1}$.\par 
 
Now suppose that we are given approximate data $\mathcal{X}_{[t,0]}(v_i)$ at the grid points \eqref{grid_points} for a band-limit $L$. If we discretize the components of the map using spherical harmonics, then we obtain a projected spherical trigonometric interpolant of the form  
\begin{equation}
\label{sph_approx_1}
\mathcal{X}_{[t,0]}  \in \mathcal{P}(\mathcal{B}_L^3 \oplus \mathcal{B}_{L-1}^3) \,.
\end{equation}In contrast, given $\mathcal{X}_{[t_h,0]}(v_i)$ and $\mathcal{X}_{[t,t_h]}(v_i)$ at the same grid points, we can perform trigonometric interpolations of both of the maps separately such that the resulting composition gives

\begin{equation}
\label{sph_approx_2}
\mathcal{X}_{[t_h,0]} \circ \mathcal{X}_{[t,t_h]} \in \mathcal{P}(\mathcal{B}_{L^2}^3 \oplus \mathcal{B}_{L^2-1}^3) \,.
\end{equation}
In turn, the approximation of the vorticity for \eqref{sph_approx_1} has an effective band-limit $L_e(t) = L \cdot L_0$ whereas for \eqref{sph_approx_2} it becomes $L_e(t) = L^2 \cdot L_0$, based on the analysis given above. Considering the composition of $N_c$ submaps, this obtains a sparse multi-scale representation of the vorticity with effective band-limit $L_e = L_0 \cdot L^{N_c}$ using only $N_c \cdot L \cdot 2L$ degrees of freedom. In comparison, an exact band-limited spherical harmonic expansion of the vorticity up to this degree would require $L^{N_c}\cdot 2L^{N_c}$ points. The spatio-temporal discretization \eqref{remapping} of the inverse map gives an exponential increase in resolution with only a linear increase in the degrees of freedom. \par

\FloatBarrier

\small
\bibliographystyle{elsarticle-num}
\bibliography{Bibliography}

\begin{thebibliography}{10}
\expandafter\ifx\csname url\endcsname\relax
  \def\url#1{\texttt{#1}}\fi
\expandafter\ifx\csname urlprefix\endcsname\relax\def\urlprefix{URL }\fi
\expandafter\ifx\csname href\endcsname\relax
  \def\href#1#2{#2} \def\path#1{#1}\fi

\bibitem{pedlosky1987geophysical}
J.~Pedlosky, et~al., Geophysical fluid dynamics, Vol. 710, Springer, 1987.

\bibitem{fjortoft1953changes}
R.~Fj{\o}rtoft, On the changes in the spectral distribution of kinetic energy
  for twodimensional, nondivergent flow, Tellus 5~(3) (1953) 225--230.

\bibitem{newton2006n}
P.~K. Newton, H.~Shokraneh, The n-vortex problem on a rotating sphere. i
  multi-frequency configurations, Proceedings of the Royal Society A:
  Mathematical, Physical and Engineering Sciences 462~(2065) (2006) 149--169.

\bibitem{zabusky1979contour}
N.~J. Zabusky, M.~Hughes, K.~Roberts, Contour dynamics for the euler equations
  in two dimensions, Journal of computational physics 30~(1) (1979) 96--106.

\bibitem{dritschel1988contour}
D.~G. Dritschel, Contour surgery: a topological reconnection scheme for
  extended integrations using contour dynamics, Journal of Computational
  Physics 77~(1) (1988) 240--266.

\bibitem{dritschel1989contour}
D.~G. Dritschel, Contour dynamics and contour surgery: numerical algorithms for
  extended, high-resolution modelling of vortex dynamics in two-dimensional,
  inviscid, incompressible flows, Computer Physics Reports 10~(3) (1989)
  77--146.

\bibitem{bosler2014lagrangian}
P.~Bosler, L.~Wang, C.~Jablonowski, R.~Krasny, A lagrangian particle/panel
  method for the barotropic vorticity equations on a rotating sphere, Fluid
  Dynamics Research 46~(3) (2014) 031406.

\bibitem{perlman1985accuracy}
M.~Perlman, On the accuracy of vortex methods, Journal of Computational Physics
  59~(2) (1985) 200--223.

\bibitem{dukowicz1987accurate}
J.~K. Dukowicz, J.~W. Kodis, Accurate conservative remapping (rezoning) for
  arbitrary lagrangian-eulerian computations, SIAM Journal on Scientific and
  Statistical Computing 8~(3) (1987) 305--321.

\bibitem{nordmark1991rezoning}
H.~O. Nordmark, Rezoning for higher order vortex methods, Journal of
  Computational Physics 97~(2) (1991) 366--397.

\bibitem{barba2005advances}
L.~Barba, A.~Leonard, C.~Allen, Advances in viscous vortex methods—meshless
  spatial adaption based on radial basis function interpolation, International
  Journal for Numerical Methods in Fluids 47~(5) (2005) 387--421.

\bibitem{magni2012accurate}
A.~Magni, G.-H. Cottet, Accurate, non-oscillatory, remeshing schemes for
  particle methods, Journal of Computational Physics 231~(1) (2012) 152--172.

\bibitem{bosler2013particle}
P.~A. Bosler, Particle methods for geophysical flow on the sphere., Ph.D.
  thesis (2013).

\bibitem{dritschel1997contour}
D.~G. Dritschel, M.~H. Ambaum, A contour-advective semi-lagrangian numerical
  algorithm for simulating fine-scale conservative dynamical fields, Quarterly
  Journal of the Royal Meteorological Society 123~(540) (1997) 1097--1130.

\bibitem{yin2021characteristic}
X.-Y. Yin, O.~Mercier, B.~Yadav, K.~Schneider, J.-C. Nave, A characteristic
  mapping method for the two-dimensional incompressible {Euler} equations,
  Journal of Computational Physics 424 (2021) 109781.

\bibitem{yin2023characteristic}
X.-Y. Yin, K.~Schneider, J.-C. Nave, A {Characteristic Mapping Method} for the
  three-dimensional incompressible {Euler} equations, Journal of Computational
  Physics (2023) 111876.

\bibitem{taylor2023projection}
S.~Taylor, J.-C. Nave, A projection-based characteristic mapping method for
  tracer transport on the sphere, Journal of Computational Physics (2023)
  111905.

\bibitem{nave2010gradient}
J.-C. Nave, R.~R. Rosales, B.~Seibold, A gradient-augmented level set method
  with an optimally local, coherent advection scheme, Journal of Computational
  Physics 229~(10) (2010) 3802--3827.

\bibitem{mercier2020characteristic}
O.~Mercier, X.-Y. Yin, J.-C. Nave, The characteristic mapping method for the
  linear advection of arbitrary sets, SIAM Journal on Scientific Computing
  42~(3) (2020) A1663--A1685.

\bibitem{arnold1966differential}
V.~Arnold, Sur la g{\'e}om{\'e}trie diff{\'e}rentielle des groupes de {Lie} de
  dimension infinie et ses applications {\`a} l'hydrodynamique des fluides
  parfaits, in: Annales de l'institut Fourier, Vol.~16, 1966, pp. 319--361.

\bibitem{marsden1983coadjoint}
J.~Marsden, A.~Weinstein, Coadjoint orbits, vortices, and clebsch variables for
  incompressible fluids, Physica D: Nonlinear Phenomena 7~(1-3) (1983)
  305--323.

\bibitem{holm1998euler}
D.~D. Holm, J.~E. Marsden, T.~S. Ratiu, The {Euler--Poincar{\'e}} equations and
  semidirect products with applications to continuum theories, Advances in
  Mathematics 137~(1) (1998) 1--81.

\bibitem{cotter2013noether}
C.~J. Cotter, D.~D. Holm, On {Noether’s} theorem for the
  {Euler--Poincar{\'e}} equation on the diffeomorphism group with advected
  quantities, Foundations of Computational Mathematics 13~(4) (2013) 457--477.

\bibitem{pavlov2011structure}
D.~Pavlov, P.~Mullen, Y.~Tong, E.~Kanso, J.~E. Marsden, M.~Desbrun,
  Structure-preserving discretization of incompressible fluids, Physica D:
  Nonlinear Phenomena 240~(6) (2011) 443--458.

\bibitem{gawlik2011geometric}
E.~S. Gawlik, P.~Mullen, D.~Pavlov, J.~E. Marsden, M.~Desbrun, Geometric,
  variational discretization of continuum theories, Physica D: Nonlinear
  Phenomena 240~(21) (2011) 1724--1760.

\bibitem{barrera1985vector}
R.~G. Barrera, G.~Estevez, J.~Giraldo, Vector spherical harmonics and their
  application to magnetostatics, European Journal of Physics 6~(4) (1985) 287.

\bibitem{grohs2019projection}
P.~Grohs, H.~Hardering, O.~Sander, M.~Sprecher, Projection-based finite
  elements for nonlinear function spaces, SIAM Journal on Numerical Analysis
  57~(1) (2019) 404--428.

\bibitem{grohs2013projection}
P.~Grohs, M.~Sprecher, Projection-based quasiinterpolation in manifolds, SAM
  Report 23 (2013).

\bibitem{gawlik2018embedding}
E.~S. Gawlik, M.~Leok, Embedding-based interpolation on the special orthogonal
  group, SIAM Journal on Scientific Computing 40~(2) (2018) A721--A746.

\bibitem{lai2007spline}
M.-J. Lai, L.~L. Schumaker, Spline functions on triangulations, Vol. 110,
  Cambridge University Press, 2007.

\bibitem{alfeld1996fitting}
P.~Alfeld, M.~Neamtu, L.~L. Schumaker, Fitting scattered data on sphere-like
  surfaces using spherical splines, Journal of Computational and Applied
  Mathematics 73~(1) (1996) 5--43.

\bibitem{alfeld1996bernstein}
P.~Alfeld, M.~Neamtu, L.~L. Schumaker, {Bernstein-B{\'e}zier} polynomials on
  spheres and sphere-like surfaces, Computer Aided Geometric Design 13~(4)
  (1996) 333--349.

\bibitem{hielscher2023approximating}
R.~Hielscher, L.~Lippert, Approximating the derivative of manifold-valued
  functions, Journal of Approximation Theory 285 (2023) 105832.

\bibitem{hairer2006structure}
E.~Hairer, C.~Lubich, G.~Wanner, Structure-preserving algorithms for ordinary
  differential equations, Geometric numerical integration 31 (2006).

\bibitem{mcewen2011novel}
J.~D. McEwen, Y.~Wiaux, A novel sampling theorem on the sphere, IEEE
  Transactions on Signal Processing 59~(12) (2011) 5876--5887.

\bibitem{boffetta2012two}
G.~Boffetta, R.~E. Ecke, et~al., Two-dimensional turbulence, Annual review of
  fluid mechanics 44~(1) (2012) 427--451.

\bibitem{arnold2021topological}
V.~I. Arnold, B.~A. Khesin, Topological methods in hydrodynamics, Vol. 125,
  Springer Nature, 2021.

\bibitem{moresi2019stripy}
L.~Moresi, B.~Mather, Stripy: A python module for (constrained) triangulation
  in cartesian coordinates and on a sphere., Journal of Open Source Software
  4~(38) (2019) 1410.

\bibitem{renka1997algorithm}
R.~J. Renka, Algorithm 772: Stripack: Delaunay triangulation and voronoi
  diagram on the surface of a sphere, ACM Transactions on Mathematical Software
  (TOMS) 23~(3) (1997) 416--434.

\bibitem{libigl}
A.~Jacobson, D.~Panozzo, et~al., {libigl}: A simple {C++} geometry processing
  library, https://libigl.github.io/ (2018).

\bibitem{haurwitz1940motion}
B.~Haurwitz, The motion of atmospheric disturbances on the spherical earth, J.
  mar. Res 3~(5) (1940) 254--267.

\bibitem{neamtan1946motion}
S.~Neamtan, The motion of harmonic waves in the atmosphere, Journal of
  Atmospheric Sciences 3~(2) (1946) 53--56.

\bibitem{lorenz1972barotropic}
E.~N. Lorenz, Barotropic instability of {Rossby} wave motion, Journal of
  Atmospheric Sciences 29~(2) (1972) 258--265.

\bibitem{tung1981barotropic}
K.~K. Tung, Barotropic instability of zonal flows, Journal of Atmospheric
  Sciences 38~(2) (1981) 308--321.

\bibitem{dritschel2015late}
D.~G. Dritschel, W.~Qi, J.~Marston, On the late-time behaviour of a bounded,
  inviscid two-dimensional flow, Journal of Fluid Mechanics 783 (2015) 1--22.

\bibitem{lindborg2022two}
E.~Lindborg, A.~Nordmark, Two-dimensional turbulence on a sphere, Journal of
  Fluid Mechanics 933 (2022).

\bibitem{kraichnan1967inertial}
R.~H. Kraichnan, Inertial ranges in two-dimensional turbulence, The Physics of
  Fluids 10~(7) (1967) 1417--1423.

\bibitem{leith1968diffusion}
C.~E. Leith, Diffusion approximation for two-dimensional turbulence, The
  Physics of Fluids 11~(3) (1968) 671--672.

\bibitem{batchelor1969computation}
G.~K. Batchelor, Computation of the energy spectrum in homogeneous
  two-dimensional turbulence, The Physics of Fluids 12~(12) (1969) II--233.

\bibitem{modin2022canonical}
K.~Modin, M.~Viviani, Canonical scale separation in two-dimensional
  incompressible hydrodynamics, Journal of Fluid Mechanics 943 (2022).

\bibitem{haynes2005stratospheric}
P.~Haynes, Stratospheric dynamics, Annu. Rev. Fluid Mech. 37 (2005) 263--293.

\bibitem{gilbarg1977elliptic}
D.~Gilbarg, N.~S. Trudinger, D.~Gilbarg, N.~Trudinger, Elliptic partial
  differential equations of second order, Vol. 224, Springer, 1977.

\end{thebibliography}

\end{document}